\documentclass{article}
\usepackage{latexsym,amsthm,amssymb,amscd,amsmath}
\newcounter{alphthm}
\newcommand{\be}{\begin{equation}}
\newcommand{\ee}{\end{equation}}
\newcommand{\ben}{\begin{enumerate}}
\newcommand{\een}{\end{enumerate}}

\def\beq{\begin{equation}}
\def\eeq{\end{equation}}

\title{Weak and strong structures and the $T_{3.5}$ property
for generalized topological spaces}
\author{E. Makai, Jr., E. Peyghan and B. Samadi}
\begin{document}


\maketitle
\begin{abstract}
We investigate weak and strong structures for 
generalized topological spaces, among others products, sums, subspaces,
quotients, and the complete lattice of generalized topologies on a given set.
Also we introduce
$T_{3.5}$ generalized topological
spaces and give a necessary and sufficient condition for a 
generalized topological space to be a $T_{3.5}$ space: they are exactly
the subspaces of powers of 
a certain natural generalized topology on $[0,1]$. For spaces with at least two
points here we can have even dense subspaces.
Also, $T_{3.5}$ generalized topological spaces are exactly
the dense subspaces of compact $T_4$ generalized topological
spaces. We show that normality is productive for
generalized topological spaces. For compact generalized topological
spaces we
prove the analogue of the Tychonoff product theorem. We prove that also
Lindel\"ofness (and $\kappa $-compactness) 
is productive for generalized topological spaces. On any ordered set we
introduce a generalized topology and determine the continuous maps between
two such generalized topological spaces: for $|X|, |Y| \ge 2$ they are the
monotonous maps continuous between the respective order topologies. We
investigate the relation of sums and subspaces of generalized topological
spaces to ways of defining generalized topological spaces.
\\\\

{\bf {Keywords}}:  Generalized topology, weak and strong structures,
products, sums, subspaces, quotients,
$T_{3.5}$, normal, compact, Lindel\"of, 
$\kappa $-compact, ordered generalized topological
spaces.\footnote{ 2010 
Mathematics subject Classification:
Primary: 54A05. Secondary: 54B30.}

\end{abstract}


\section{Introduction}

In this paper we do not require acquaintance with the terminology of category
theory, although we use some of its concepts. These will be explained in the
respective places.

In \S 2 we collect material needed later in our paper, and give the necessary
definitions.

In \S 3 we investigate weak and strong structures for generalized topological
spaces (GTS's), in
particular, products (different from Cs\'asz\'ar's products), sums, subspaces,
quotients, and the complete lattice of generalized topologies (GT's) 
on a set $X$. 
Our definition of the product is the categorical definition. 
It will turn out that GTS's form a topological category over sets (its
definition cf. in Theorem 3.1.) This is a slight variant of \cite{KM}, Theorem
4.8.

In \S 4 we will investigate
productivity of certain topological properties with respect to our definition
of product. These include the natural analogues for GTS's
of the separation axioms $T_0$, $T_1$, $T_2$ ,$T_3$, but also that
of $T_4$.
For compact GTS's
there holds the
analogue of the Tychonoff product theorem. However, also Lindel\"of property
and $\kappa $-compactness are productive for GT's.
We will define $T_{3.5}$
GTS's that have an analogous relation to the GT on $[0,1]$
having a base $\{ [0,x),\,\,(y, 1] \mid x,y \in [0,1] \} $ as 
$T_{3.5}$ topological spaces (TS's) 
have to the usual topology on $[0,1]$: a GTS is $T_{3.5}$ if and
only if it is a subspace of some power of the GTS $[0,1]$ if and only if it is
a subspace of a normal $T_4$ GTS. For ordered spaces $(X, \le )$
there is a natural GT on $X$,
and the continuous functions between two such spaces $X,Y$, for $|X|, |Y| \ge
2$ are exactly the
monotonous maps continuous in the respective order topologies.

In \S 5 we will
investigate the relation of generating GT's by a monotonous map $\gamma :
P(X) \to P(X)$, and by an enlargement $k: \mu \to P(X)$, to subspaces and sums 
of GTS's.

 
\section{Preliminaries}


{\bf{1.}}
The concept of generalized topology dates back to antiquity, then
called ``closure operator'' (which could have still some additional
properties, like, e.g., idempotence). A large number of such additional
properties of closure operators and their
interrelations are discussed in the monographs \cite{DT} and \cite{Ca}.
Early examples are the linear spans of a subset of a vector
space, or more generally, subalgebras generated by 
subsets of some algebraic structure,
like groups,  semigroups, etc.
For history and many properties of such closure operators 
cf. the papers \cite{DG} and \cite{DGT} from 1987 and 1989, and particularly
the monograph of D. Dikranjan and W. Tholen \cite{DT} from 1995, and 
the more recent monograph of G. Castellini \cite{Ca} from 2003. 
Also cf.
the monograph of E. \v Cech, Z. Frol\'\i k, M. Kat\v etov \cite{CFK}, from
1966, but that deals only with one type of closure spaces, 
called \v Cech-closure spaces, or pretopologies (definition cf. later).

Let $X$ be a set and $P(X)$ its power set.
\cite{DT}, pp. (xiii) and 147, defined a {\it{closure operator}} 
$c:P(X) \to P(X)$ as follows. It
should be {\it{increasing}} (called there {\it{extensive}}) i.e.,
$A \subset cA$ and {\it{monotonous}} 
i.e., $A \subset B \Longrightarrow cA \subset cB$. A {\it{closure space}},
also written as {\it{CS}}, is a pair $(X,c)$, where $X$ is a set and $c:P(X) \to
P(X)$ is a closure operator.
\cite{DT}, p. 147 
also investigated {\it{continuous maps between closure spaces}} 
$f:(X,c) \to (Y,d)$, i.e., maps $X
\to Y$, satisfying 
\begin{equation}\label{2.1}
A \subset X \Longrightarrow fcA \subset dfA \Longleftrightarrow cA \subset
f^{-1}dfA\,,
\end{equation}
or, equivalently,
\begin{equation}\label{2.2} 
B \subset Y \Longrightarrow c f^{-1}B \subset f^{-1}dB \Longleftrightarrow
fcf^{-1}B \subset dB
\end{equation}
(cf. \cite{DT}, p. 25 and \cite{Ca}, p. 42, Proposition 4.2). 
All closure spaces and all continuous maps between them form a (so
called) {\it{category}}, denoted by {\bf{CS}}. 
(Actually the setting of \cite{DT} and \cite{Ca}
was more general: a category $\cal
X$, with a distinguished class $\cal M$ of subobjects, and the closure
operator mapped any distinguished subobject of any object $X$ of $\cal X$
to some distinguished subobject of the same object $X$. 
Additionally, all morphisms $f:X \to Y$ were required to be continuous from the
closure operator on $X$ to the closure operator 
on $Y$. E.g., for topological groups, each
morphism $f$ carries the closure of any subgroup $X_0$ of $X$ 
into the closure of the subgroup $f(X_0)$ of $Y$. Here $\cal X$ and
$\cal M$ had to satisfy some natural hypotheses, which hold in our
cases. However, their main topic is not a generalization of the investigation of
generalized topological spaces. Namely, for ${\cal X} = {\bold{Set}}$ and $\cal
M$ being all monomorphisms in ${\bf{Set}}$
their resulting category is just ${\bf{Set}}$. If we let
${\cal X} = {\bf{GenTop}}$
and $\cal M$ all monomorphisms in ${\bf{GenTop}}$, the category
${\bf{GenTop}}$ is already contained in the hypotheses, so this is no
definition of ${\bf{GenTop}}$.
A reader not interested in category theory may just skip this point.)
Initial, i.e., weak
and final, i.e., strong structures for supratopological spaces --- which are
closely related to generalized topological spaces, cf. below --- are
proved to exist and are investigated in \cite{KM}.

\cite{CFK} required that a closure operator $c:P(X) \to P(X)$ should be
increasing, 
and {\it{preserve finite unions}}, also called {\it{finitely additive}},
i.e., $c \emptyset 
= \emptyset $ (in \cite{DT} p. xiii {\it{groundedness}})
and $A, B \subset X \Longrightarrow c(A \cup B)=(cA) \cup (cB)$ (in \cite{DT}
p. xiii and in \cite{Ca}, p. 65, Definition 6.1
{\it{additivity}}).
Such an operator $c$ is called a {\it{\v Cech-closure}} and the pair $(X,c)$ a
{\it{\v Cech-closure space}}, or more recently a {\it{pretopology}} and a
{\it{pretopological space}}.
The pretopological spaces with the corresponding continuous maps were
investigated in great detail in \cite{CFK}. They form the (so called) 
{\it{category}} {\bf{PrTop}}. In particular, initial, i.e., weak
and final, i.e., strong structures for pretopological spaces are
proved to exist and are investigated in detail in \cite{CFK}, 
\S 32 and in \S 33.
In our paper pretopological spaces will not be investigated. We have to
remark that also in \cite{DT} most of the concrete examples in topology
were connected
with pretopological spaces, while in \cite{Ca} such examples are rare --- just
pretopological spaces are defined in p. 91, Example 7.12 ---
and GTS's and closure spaces were not
systematically investigated from the topological point of view in \cite{DT}
and \cite{Ca}.

{\bf{2.}}
In topology, generalized topologies $(X, \mu )$ formally seem (almost)
to have been defined by A. S. Mashhour,
A. A. Allam, F. S. Mahmoud, F. H. Khedr \cite{MAMKh}, in 1983,
under the name of ``{\it{supratopological spaces}}'', where however $X$ open was
required (a strong generalized topology). This concept and terminology persist
till now. Categorical topologists investigate them, as one of the many types
of structures in topology (the most well-known of these are beside topological
spaces the uniform
spaces), and investigate the relationships of these different types of 
structures in topology. However, unfortunately the terminologies collide:
categorical topologists
used to call supratopological spaces also as closure spaces (cf. e.g.,
\cite{KM}), which is in
conflict with the usage of the monograph \cite{DT}. We will use the term
{\it{strong generalized topological space}}.
What is closure space in \cite{DT}, yet satisfying the extra condition that the
closure of the empty set is the empty set, is called by categorical
topologists a {\it{neighbourhood
space}}. It is given by a system of {\it{neighbourhoods}} ${\cal{N}}(x)$
for each point of a set $X$, such that $N \in {\cal{N}}(x) \Longrightarrow
x \in N$, $N' \supset N \in {\cal{N}}(x) \Longrightarrow N' \in {\cal{N}}(x)$,
and $X \in {\cal{N}}(x)$.

Then \'A. Cs\'asz\'ar \cite{Cs1} in 2002
introduced generalized topological spaces, which differ from supratopological
spaces just by omitting the requirement 
of openness of $X$ from their definition.
His motivation was the previous investigation of a number of generalizations
of open sets in topological spaces, like semiopen sets (\cite{Le}, 1963), 
$\alpha $-open sets (\cite{Nj}, 1965), preopen sets (\cite{MAE}, 1982), 
$\beta $-open sets (\cite{AEM}, 1983) defined by 
$A \subset {\text{cl\,int}}\, A$, 
$A \subset {\text{int\,cl\,int}}\,A$,
$A \subset {\text{int\,cl}}\,A$, 
$A \subset {\text{cl\,int\,cl}}\,A$, 
for $A$ a
subset of a topological space $X$, respectively. 
An extensive literature cf. in \cite{Cs1} from 2002. 
These definitions led \'A. Cs\'asz\'ar \cite{Cs0.5} to introduce their
common generalization, the so called $\gamma $-open sets, where $\gamma : P(X)
\to P(X)$ is an arbitrary monotonous map, via the property $A \subset \gamma A$.
The concept of $\gamma $-open sets already includes all generalized
topologies (for suitable $\gamma $, namely for $\gamma $
the interior operator of the generalized topology). In \cite{Cs1}
\'A. Cs\'asz\'ar made a further step: he considered the system of $\gamma
$-open sets, which is always a generalized topology, and disregarded from
which $\gamma $ was it derived. Thus he \cite{Cs1}
arrived to the concept of generalized topologies, and
began their systematic topological investigation.
The paper \cite{Cs1} was the basis for at least 400 subsequent
papers in this subject (by MathSciNet).
This has been one of the important developments of general topology in the
recent years. 

We note that beginning with a topology,
the first four above given generalizations of open
sets form only generalized topologies, except for $\alpha $-openness.
Moreover, the first four above
types of generalized open sets can be introduced also in generalized
topological spaces, cf. \cite{Cs1.6}.

We remark that the difference between supratopological spaces and generalized
topological spaces is minor. Many proofs for supratopologies carry over to
generalized topologies, sometimes with some notational complications. But of
course, there are also differences between them.
 
\'A. Cs\'asz\'ar \cite{Cs4} in 2008 introduced {\it{generalized neighbourhood
systems}}, which is a generalization of the above mentioned neighbourhood
spaces, by omitting the condition that $X \in {\cal{N}}(x)$.
This concept is equivalent to that of the closure spaces. In fact, from a
closure operator $c:P(X) \to P(X)$ one derives ${\cal{N}}(x)$ by $N \in
{\cal{N}}(x) \Longleftrightarrow x \not\in c(X \setminus N)$, and the same
formula derives $c:P(X) \to P(X)$ from $\langle {\cal{N}}(x)
\mid x \in X \rangle $. Continuity can be rewritten as follows: $f: (X, \langle 
{\cal{N}}(x) \mid x \in X \rangle ) \to (Y,\langle {\cal{M}}(y) \mid y \in Y
\rangle )$ satisfies $x \in X \Longrightarrow f^{-1}{\cal{M}}(fx) \subset 
{\cal{N}}(x)$.

{\bf{3.}}
Let $X$ be a set and $\mu\subset P(X)$. (We observe that 
some authors require still $X \ne \emptyset $.  
However then e.g. intersections of subspaces are not subspaces, the empty sum
does not exist, etc.,
so we must allow $X = \emptyset $.) 
Then $\mu$
is called a {\it{generalized topology}}, briefly {\it{GT}} on $X$ if 
$\emptyset \in\mu$ and any
union of elements of $\mu $ belongs to $\mu $. A set $X$ with a GT $\mu $
is said to be a
{\it{generalized topological space $(X, \mu )$}}, briefly {\it{GTS}}. 
The elements of $\mu$ are called
$\mu$-{\it{open sets}}, and their complements are called
{\it{$ \mu $-closed}}. 
We say that $\mu$ is {\it{strong}} if $X\in\mu$. 
A {\it{base of a GTS}} $(X, \mu )$ is a subset $\beta
$ of $\mu $ such that each $M \in \mu $ is a union of a subfamily (possibly
empty) of $\beta $, cf. \cite{Cs1.7}.

For $A\subset X$, we
denote by $c_\mu(A)$ the intersection of all $\mu$-closed sets containing $A$
and by $i_\mu(A)$ the union of all $\mu$-open sets contained in
$A$. Then the map $c_{\mu }: P(X) \to P(X)$ is increasing, 
monotonous and {\it{idempotent}} (i.e., $c^2=c)$. If some
$c : P(X) \to P(X)$ has these properties, then it defines a GT via $\mu _c :=
\{ X \setminus c(A) \mid A \subset X \} $. The description of GT's by open
sets, or by the closure operator are equivalent: $\mu $ is sent to $c_{\mu }$,
and $c$ to $\mu _c$, and these maps define bijections inverse to each other. 
For GT's we will use the notations $\mu, \nu, \varrho $ for the set of all
open sets, and the notations
$c,d,e$ for the associated closure operators.
(The description by closed sets is clearly equivalent to
the description by open sets, so we will not consider it in this paper.)

For maps, $f:(X, \mu ) \to (Y, \nu )$ or $f:(X,c) \to (Y,d)$
is {\it{continuous}} if $f^{-1}(\nu )\subset \mu $, or in terms of closure
operators if (\ref{2.1})
or, equivalently, (\ref{2.2}) holds.
We will write also that $f$ is $(\mu , \nu )$-{\it{continuous}}, or
$(c,d)$-{\it{continuous}}.
Identifying the $\mu $'s and $c$'s on a set $X$ 
via the above bijections $\mu \mapsto c_{\mu }$ and $c \mapsto \mu _c$, 
these concepts become
equivalent. The GTS's, with the continuous maps between any two of them form a
(so called) {\it{category}}, denoted by {\bf{GenTop}}. If $f:(X, \mu ) \to (Y,
\nu )$ is continuous, and has a continuous inverse $g:(Y, \nu ) \to (X, \mu
)$ (or we may use $(X,c),(Y,d)$), then it is called a {\it{homeomorphism}}.

For ${\cal{A}} \subset P(X)$ and $X_0 \subset X$ we write ${\cal{A}}|X_0 :=
\{ A \cap X_0 \mid A \in {\cal{A}} \} $ ({\it{trace of ${\cal{A}}$ on $X_0$}}).

{\bf{4.}}
For concepts of category theory, we refer to \cite{HS}, \cite{AHS}.
However, 
in this paper we do not want to suppose acquaintance with category theory.
An exception is when we speak about limits or colimits, but then
the reader may restrict himself to their special cases products or sums.
However, because of 
this we have to recall the general concepts of products and sums
in categories. 

We begin with a notation.
Recall that $\{ ( \cdot )_{\alpha } \mid \alpha \in J \} $ 
is the family (set, or class) 
of all $( \cdot )_{\alpha}$'s, for $\alpha \in J$. Here
multiple occurrence of the same $( \cdot )_{\alpha }$ amounts to the same as
if it occurred only once.
If we write $\langle ( \cdot )_{\alpha } \mid \alpha \in J \rangle $, this
means the {\it{indexed family}} 
of all $( \cdot )_{\alpha}$'s, for $\alpha \in J$.
That is, $\langle ( \cdot )_{\alpha } \mid \alpha \in J \rangle $ is a
function from $J$, whose values may coincide for different $\alpha $'s.
When the indexed family is a set (i.e., $J$ is a set), 
we write {\it{indexed set}}.

In a category, like e.g. that of all sets (as objects)
and all functions
between them (as morphisms), or all generalized topological spaces (as
objects) and all continuous maps
between them (as morphisms), etc., 
one defines {\it{products}} and {\it{sums}} 
in the following way. 

{\it{For an indexed set of objects $\langle
X_{\alpha } \mid \alpha \in J \rangle $ 
their product}} $\prod _{\alpha \in J} X_{\alpha
}$ is the up to isomorphism unique object, for which there exist so called
{\it{projections}} $ \pi _{\alpha }: \prod _{\alpha \in J} X_{\alpha } \to 
X_{\alpha }$, which have the following universality property. For any
morphisms $\langle f_{\alpha }: Y \to X_{\alpha } \mid \alpha \in J\rangle $ 
there exists a unique morphism $g:Y \to \prod _{\alpha \in J} X_{\alpha }$
such that for each $\alpha \in J$ we have $f_{\alpha }=\pi _{\alpha} g$.
The {\it{underlying set}} of $(X, \mu )$ (or of $(X,c)$) is $X$.
(More details cf. in \S 3.) 
In our categories the product of an indexed set of objects may be supposed to
have as underlying set the product of the underlying sets. 
We will actually suppose this. 

{\it{For an indexed set of objects $\langle
X_{\alpha } \mid \alpha \in J \rangle $ 
their sum\,}} $\coprod _{\alpha \in J} X_{\alpha
}$ is the up to isomorphism unique object, for which there exist so called
{\it{injections}} $ \iota _{\alpha }: X_{\alpha } \to 
\coprod _{\alpha \in J} X_{\alpha }$, 
which have the following universality property. For any
morphisms $\langle f_{\alpha }: X_{\alpha } \to Y \mid \alpha \in J \rangle $ 
there exists a unique morphism $g:\coprod _{\alpha \in J} X_{\alpha } \to Y$
such that for each $\alpha \in J$ we have $f_{\alpha }=g \iota _{\alpha} $.
In our categories the sum of an indexed set of objects may be supposed to
have as underlying set the sum (i.e., disjoint union) 
of the underlying sets. We will actually
suppose this. Moreover, we may identify $X_{\alpha } $ and $\iota _{\alpha
}X_{\alpha }$ via $\iota _{\alpha }$, 
and then we may consider the underlying set of $\coprod
_{\alpha \in J} X_{\alpha }$ as the disjoint union of the underlying sets of
the $X_{\alpha }$'s, that we will do also.

Analogously, if we have an injection $X_0 \to X$, then we may 
consider this as an inclusion of a subset, that we will do as well.

The {\it{empty product}} (or {\it{$0$'th power of a space}}) 
is by this definition that up to isomorphism unique
object $X_{\text{fin}}$ ({\it{final object}})
for which for any object $X$ there is exactly one
morphism $X \to X_{\text{fin}}$. In {\bf{GenTop}} this
is $( X, \{ \emptyset \} )$, where $|X|=1$. 
Similarly, the {\it{empty sum}} is that object
$X_{\text{init}}$ ({\it{initial object}})
for which for any object $X$ there is exactly one
morphism $X_{\text{init}} \to X$. In {\bf{GenTop}} this
is $( \emptyset , \{ \emptyset \} )$.

{\bf{5.}}
A way to produce GT's is given by the following (\cite{Cs1}).
We call $\gamma:P(X) \rightarrow P(X)$ monotonous as in \S 2, {\bf{1}}
(and write $\gamma
A$ for $\gamma(A)$), and denote by $\Gamma (X)$ the family of all such
mappings.
A set $A \subset X$ is said to be $\gamma $-{\it{open}} 
if $A \subset \gamma
A$. The $\gamma$-open sets constitute a GT on $X$ (cf. \cite{Cs0.5}, 1.1), 
which we denote by $\mu ( \gamma )$. Actually, all
GT's on a given set $X$ can be obtained in this way (see Lemma 1.1 of
\cite{Cs1}).

Another way to produce GT's is given by the following (see
\cite{Cs3}). A mapping $k:\mu \rightarrow P(X)$ is said to be an 
{\it{enlargement on
$(X, \mu )$}} if $M\subset kM$, whenever $M\in \mu$. A subset $A\subset X$ is
{\it{$\kappa (\mu ,k)$-open}} iff $x\in A$ implies the the existence of
 a $\mu$-open set $M$ such that $x\in M$ and $kM\subset A$. 
Cs\'{a}sz\'{a}r
 in \cite{Cs3} proved that the collection $\kappa (\mu ,k)$ 
of all $\kappa (\mu ,k)$-open sets
 is a GT on $X$ that is {\it{coarser than}} $\mu$ 
(i.e. $\kappa (\mu ,k) \subset \mu$) whenever $\mu$ is a GT on $X$.
Some further aspects of enlargements are investigated in Y. K. Kim,
W. K. Min \cite{KimMin}.

\cite{Cs5} defined a sort of product of GTS's and obtained some of its basic
properties. One can find more results related to this concept in \cite{M},
\cite{PC} and \cite{Sh}. Its definition cf. in \S 3 of this paper, Definition
3.11. We will call this the {\it{Cs\'asz\'ar product of GTS's}}, but this
will not be investigated in our paper.

{\bf{6.}}
There are some papers related to separation axioms on GTS's such as
\cite{Cs1.5}, \cite{M}, \cite{Cs2} and \cite{XY}. 
In particular, $T_0$, $T_1$, $T_2$, {\it{regularity}}, $T_3$ 
(i.e., regular $T_1$, or
equivalently regular $T_0$, like for topologies),
{\it{normality and}} $T_4$ (i.e., normal $T_1$) are defined word for word as
for topological spaces. (Observe that if $X$ is a $T_1$ --- e.g., $T_2$ --- 
GTS with $|X| \ge 2$
or $X = \emptyset $ then $X$ is strong. For $|X|=1$ there are two GT's on $X$:
$(X, \{ \emptyset \} )$ and $(X, \{ \emptyset , X \} )$. Both are $T_2$ hence
$T_1$, and the
first one is not strong, the second one is strong. For normality the situation
is converse: a not strong GTS is vacuously 
normal, since there are no two disjoint closed subsets (empty or non-empty). 
For regularity we
have: a GTS of the form
$(X, \{ \emptyset \} )$ is vacuously regular --- but regularity of a
GTS $(X, \mu )$ with $\mu \ne \{ \emptyset \} $ implies strongness of 
$(X, \mu )$. Normal $T_0$ does not imply $T_1$, already for topologies, e.g.
for ${\mathbb{R}}$ with open base $\{ ( -\infty , r) \mid r \in {\mathbb{R}}
\} $.)
Also, \cite{Cs2} studied normal GTS's and exhibited
a suitable form of Urysohn lemma (\cite{Cs2}, Theorem 3.3)
by defining a suitable GT on $[0,1]$: it has as a base $\{ [0,x), (y,1] \mid
x,y \in [0,1] \} $.

{\bf{7.}}
Let $X$ be a set. We say that ${\cal{A}} \subset P(X)$ is a {\it{stack}} 
(also called {\it{ascending)}} if $A \in {\cal{A}}$ and $A \subset B \subset X$
imply $B \in {\cal{A}}$ (cf. \cite{KM}). For ${\cal{B}} \subset P(X)$ we write 
Stack\,${\cal{B}}:= \{ C \subset X \mid \exists B \in {\cal{B}} {\text{ such
that }} C \supset B \} $.
This is called the {\it{stack generated by}} ${\cal{B}}$ (or the {\it{ascending
hull of}} ${\cal{B}}$).

We recall that \cite{Cs4} investigated {\it{generalized neighbourhood
systems}} (GNS's), 
i.e., functions $\psi :X \to P \left( P(X) \right) $ with the property that
$x \in X$ and $V \in \psi (x)$ imply $x \in V$. The filter property is not
required, and also $\psi (x) = \emptyset $ is allowed. Moreover, not even the
stack property
is required. However, for what \cite{Cs4} uses these
GNS's, remains invariant if we take the
generated stack Stack\,$\psi (x):=($Stack\,$\psi )(x):=$Stack\,$(\psi x)$
rather than $\psi (x)$. 
That is, we may suppose that each $\psi
(x)$ is a stack. These functions $\psi ( \cdot )$ are equivalent to the closure
operators defined above, cf. the end of \S 2, {\bf{2}}.

\cite{Cs4} defined for such a
$\psi $ a GT $\mu _ {\psi }:= \{ M \subset X \mid x \in M \Longrightarrow 
\exists V \in \psi (x)\,\,\,\,V \subset M\} $, called the {\it{GT generated 
by}} $\psi $. By \cite{Cs4}, Lemma 2.2 here for ${\text{Stack}}\,\psi $ 
we have $\mu _{\psi } = \mu _{{\text{Stack}}\,\psi }$, so here 
we may assume the stack property of $\psi $. \cite{Cs4}, p. 396
established that each GT can be generated by at least one GNS, that can be
supposed by \cite{Cs4}, Lemma 2.2 to consist of stacks. \cite{Cs4}   
Example 2.1 and p. 397 showed that several different 
GNS's $\psi =$ Stack\,$\psi $ 
can generate the same $\mu _{\psi }$. Conversely, {\it{a GT $\mu $
generates a GNS}}, by the formula $\psi _{\mu }(x) 
:= {\text{Stack}}\,\{ M \in \mu \mid x \in M \}$, \cite{Cs4}, 
proof of Lemma 1.3.
For two GNS's there can be defined the
continuous maps: $f: (X, \psi ) \to (Y, \psi ')$ {\it{is continuous}} if and
only if $x \in X \Rightarrow f^{-1}
\left( \psi ' \left( f(x) \right) \right) \subset \psi (x)$. 
These continuous maps
remain continuous if we replace the GNS's by the GT's generated by them, cf. 
\cite{Cs4}, Proposition 2.1. However, the converse is not true:
a map between GNS's which is continuous between the generated GT's is not
necessarily continuous between the GNS's. 
Cf. \cite{Cs1}, Example 2.2, where the generated GT's are even equal, so
different GNS's may generate the same GT. 
These discrepancies between GT's and GNS's are the difference between the 
categories of GTS's and CS's.

W. K. Min \cite{Mi} investigated the relationship of GNS's and GT's further.


\section{Topologicity of {\bf{GenTop}} over {\bf{Set}}}

Let {\bf{GenTop}} be the {\it{category}} of all GTS's (called
{\it{objects}}) and all continuous maps between them
(called {\it{morphisms}}). Similarly, {\bf{Set}} is the
category of all sets (as objects) and all functions between sets (as morphisms).
As well known, there are several ways to define the
category {\bf{GenTop}}. E.g., with generalized open sets, i.e.,
objects are pairs $(X, \mu )$, with $\{ \emptyset \} \subset \mu 
\subset P(X)$, where $\mu $ is closed under arbitrary unions, 
and morphisms $f:(X, \mu ) \to (Y, \nu )$ characterized by
$f^{-1} \nu \subset \mu $. Or with closure operators $c:P(X) \to P(X)$,
that are increasing, monotonous and idempotent.
Then, as for topological spaces, $f:(X,c) \to (Y,d)$ is a
morphism iff (\ref{2.1}) holds, 
or, equivalently, iff (\ref{2.2}) holds
(cf. beside \cite{DT}, p. 25 and \cite{Ca}, p. 42, Proposition 4.2
also \cite{Cs1}, \cite{Sa12}).

A {\it{source}}, or {\it{sink}} in a category ${\cal{C}}$ is an indexed family
(set or class) 
of morphisms with common domain, or codomain, i.e., $\langle f_{\alpha
}:X \to Y_{\alpha } \mid \alpha \in J \rangle $, or $\langle g_{\alpha }:
X_{\alpha } \to Y \mid \alpha \in J \rangle $.

If we have an {\it{indexed family of mappings}} 
$\langle f_{\alpha }: P(X) \to P(Y)
\mid \alpha \in J \rangle $, for some sets $X$ and $Y$, then
their {\it{union and intersection}} are defined pointwise: 
\begin{equation}\label{3.1}
\begin{array}{l}
(\cup _{\alpha \in
J} f_{\alpha })(A) := \cup _{\alpha \in J} (f_{\alpha }A), {\text{ and }} \\
(\cap _{\alpha \in
J} f_{\alpha })(A) := \cap _{\alpha \in J} (f_{\alpha }A), {\text{ for }} 
A \subset X \,.
\end{array}
\end{equation}

The following theorem has to be preceded by some definitions.
The {\it{underlying set functor}} $U:$ {\bf{GenTop}} $ \to $ {\bf{Set}} 
maps the generalized topological space $(X,
\mu )$ (or $(X,c)$) to the set $X$ and the continuous map 
$f:(X, \mu ) \to (Y, \nu )$ (or $(X,c) \to (Y,d)$)
to the function $f:X \to Y$. 
$U$ is called
{\it{faithful}} if $f,g:(X, \mu ) \to (Y, \nu )$ being different implies that
also $Uf,Ug:X \to Y$ are different. This is evident for {\bf{GenTop}}. 
$U$ is called
{\it{amnestic}} if the following holds. 
If the identity map on $X$ is (underlies) a continuous map $(X, \mu) \to (X,
\nu ) $ and also is (underlies) a continuous map $(X, \nu) \to (X,
\mu ) $, then $\mu = \nu $. This is also evident for {\bf{GenTop}}.
$U$ is
{\it{fibre-small}}, if any set $X$ is the underlying set of only set many
GTS's. This is also evident, since 
the cardinality of GT's
on a set $X$ is at most $\exp (\exp|X|)$). 
                                  
The definition of {\it{initial lifts}}, also called {\it{weak structures}} 
is the following. 
If for a set $X$ there is an indexed class of
morphisms $ \langle \varphi _{\alpha }:
X \to U(Y_{\alpha }, \nu _{\alpha })= Y_{\alpha } \mid \alpha \in J \rangle $
in {\bf{Set}} 
(i.e., a source in {\bf{Set}}),
then there is a (by the way, unique) so called 
{\it{initial, or weak structure}} 
$(X, \mu )$ on $X$, such that all $\varphi _{\alpha }$
underlie 
morphisms $f_{\alpha }:(X, \mu ) \to (Y_{\alpha }, \nu _{\alpha })$,
and this source $ \langle f_{\alpha } \mid \alpha \in J \rangle $
in {\bf{GenTop}}
has the following universality property.
If for any $(Z, \lambda )$ there are morphisms
$\langle g_{\alpha } : (Z, \varrho ) \to (Y_{\alpha }, \nu _{\alpha }) \mid
\alpha \in J \rangle $ (another source in {\bf{GenTop}})
such that 
$$
Ug_{\alpha }=\varphi _{\alpha } h {\text{ for all }} \alpha \in J
$$ 
for some $h:UZ \to X$, then there
exists an $h' : (Z, \varrho ) \to (X, \mu )$, such that 
$$
h=Uh'
{\text{ and for each }} \alpha \in J  {\text{ we have }}
g_{\alpha } = f _{\alpha } h'.
$$

If we reverse in this definition the direction of the maps (i.e., $\to $ is
replaced by $\leftarrow $ and vice versa), 
we obtain the definition of {\it{final lifts}}, also called 
{\it{strong structures}}.
In details, this is the following.
If for a set $X$ there is an indexed class of
morphisms $ \langle \psi _{\alpha }:
 U(Y_{\alpha }, \nu _{\alpha })= Y_{\alpha } \to X \mid \alpha \in J
\rangle $ in {\bf{Set}} (i.e., a sink in {\bf{Set}}),
then there is a (by the way, unique) so called 
{\it{final structure}} $(X, \mu )$ on $X$, such that all $\psi _{\alpha }$
underlie 
morphisms $f_{\alpha }: (Y_{\alpha }, \nu _{\alpha }) \to 
(X, \mu )$,
and this sink $ \langle f_{\alpha } \mid \alpha \in J \rangle $
in {\bf{GenTop}}
has the following universality property.
If for any $(Z, \varrho )$ there are morphisms
$\langle g_{\alpha } : (Y_{\alpha }, \nu _{\alpha }) \to (Z, \varrho ) \mid
\alpha \in J \rangle $
(another sink in {\bf{GenTop}}) such that 
$$
Ug_{\alpha }=h \psi _{\alpha } {\text{ for all }} \alpha \in J
$$ 
for some $h:X \to UZ$, then there
exists a $h' : (X, \mu ) \to (Z, \varrho )$, such that 
$$
h=Uh'
{\text{ and for each }} \alpha \in J  {\text{ we have }}
g_{\alpha } = h' f _{\alpha } \,.
$$

The existence of all initial (weak)
structures is equivalent to the existence of all
final (strong) structures (faithfulness and fibre-smallness supposed)
\cite{AHS}, Theorem 21.9, \cite{Pr}, \cite{Enc}.

The initial lift of (i.e., weak structure for)
the empty source for a set $X$ is the indiscrete GT, i.e.,
$(X, \{ \emptyset \} )$, i.e., $(X,c)$ with 
$\forall A \subset X \,\,\,\,cA=X$.
The final lift of (i.e., strong structure for) the empty sink for a set
$X$ is the discrete GT, i.e.,
$\left( X, P(X) \right) $, i.e., $(X,c)$ with
$\forall A \subset X \,\,\,\,cA=A$.
(Therefore in our Propositions 3.2, 3.3, 3.4, 3.6 we might investigate
only non-empty sources and sinks.)

For {\it{topological categories over}} {\bf{Set}} we refer to J. Ad\'amek,
H. Herrlich, G. E. Strecker, Abstract and concrete categories: the joy of
cats, \cite{AHS}, Ch. 21, 
G. Preuss, Theory of Topological Structures \cite{Pr}, 
or for a synopsis Encyclopedia of Math. Vol. 9 \cite{Enc} 
pp. 201-202. Cf. also the text of Theorem 3.1 for the definition.


\vskip.2cm

{\bf{Theorem 3.1}} (For supratopologies cf. \cite{KM}, Theorem 4.8.)
{\it{The category {\bf{GenTop}}, with its underlying set
functor $U: {\text{\bf{GenTop}}} \to {\text{\bf{Set}}}$ 
is a topological category over {\bf{Set}}. That is, $U$ is faithful,
amnestic, 
fibre-small, and there exist all initial lifts (i.e., weak structures) or
equivalently there exist all final lifts (i.e., strong structures). 
Hence in {\bf{GenTop}}
there exist both limits and colimits of all diagrams, which can be
obtained from the respective underlying diagrams in {\bf{Set}} 
by initial/final lifts.}}
  
\vskip.2cm


Above we already observed faithfulness, amnesticity 
and fibre-smallness for $U$. Existence
of all limits (e.g., products) and all colimits (e.g., sums) and the way of
obtaining them
hold in any topological category over
{\bf{Set}} \cite{AHS}, Proposition 21.15, \cite{Enc}.
So only the weak and strong structures need be given.


We give the simple proof of Theorem 3.1, even in
several forms. We explicitly give all initial and all final lifts, i.e., weak
and strong structures, both for
the open sets and the closure operator definition. 


\vskip.2cm

{\bf{Proposition 3.2}} (For supratopologies cf. \cite{KM}, Proposition 4.1 and
Theorem 4.4.)
{\it{Let $X$ be a set. Let us have a source 
$\langle \varphi _{\alpha } :
X \to U(Y_{\alpha }, \nu _{\alpha })= Y_{\alpha } \mid \alpha \in J \rangle $ 
in {\bf{Set}}. 
Then its initial lift (i.e., weak structure) in {\bf{GenTop}} is 
$$
(X, \mu ) := (X, \{ \cup _{\alpha \in J} 
\varphi _{\alpha } ^{-1} (M_{\alpha }) \mid 
M_{\alpha } \in \nu _{\alpha } \} ) .
$$}}


\begin{proof}
Clearly $(X, \mu )$ is a GTS, and $\varphi _{\alpha }$ becomes (underlies) a
continuous map $f_{\alpha }:(X, \mu ) \to (Y_{\alpha }, \nu _{\alpha })$ 
in {\bf{GenTop}}, for each $\alpha \in J$. 

We turn to show the universality property. Let us have 
for some $(Z, \lambda )$  morphisms
$\langle g_{\alpha } : (Z, \varrho ) \to (Y_{\alpha }, \nu _{\alpha }) \mid
\alpha \in J \rangle $ (a source in {\bf{GenTop}})
such that 
$$
Ug_{\alpha }=\varphi _{\alpha } h {\text{ for all }} \alpha \in J
$$ 
for some $h:UZ \to X$. Then 
$$
\varrho \supset g_{\alpha } ^{-1} ( \nu _{\alpha }) = (U g_{\alpha }) ^{-1}
( \nu _{\alpha }) = (\varphi _{\alpha } h) ^{-1} ( \nu _{\alpha }) = 
h ^{-1} \varphi _{\alpha }  ^{-1} ( \nu _{\alpha })\,.
$$
Hence 
$$
\varrho \supset \cup _{\alpha  \in J} 
h ^{-1} \varphi _{\alpha }  ^{-1} ( \nu _{\alpha }) = 
h ^{-1} \cup _{\alpha  \in J} \varphi _{\alpha } ^{-1} ( \nu _{\alpha })\,.
$$
Thus $h ^{-1}$ maps $\mu $ into $\varrho $, hence $h = U h'$ for a
continuous map $h': (Z, \varrho ) \to (X, \mu )$.
\end{proof}


{\bf{Proposition 3.3.}} {\it{Let $X$ be a set. Let us have a sink $\langle \psi
_{\alpha }: U(Y_{\alpha }, \nu _{\alpha })= Y_{\alpha } \to X \mid \alpha \in
J \rangle $ in
{\bf{Set}}. Then its final lift (strong structure) in {\bf{GenTop}} is 
$$
(X, \mu ) := (X, \{ M \subset X
\mid \forall \alpha \in J \,\,\,\psi _{\alpha }^{-1} (M) 
\in \nu _{\alpha } \} ).
$$}}

 
\begin{proof}
Clearly $(X, \mu )$ is a GTS, and $\varphi _{\alpha }$ becomes (underlies) a
continuous map $f_{\alpha }: (Y_{\alpha }, \nu _{\alpha }) \to
(X, \mu )$ in {\bf{GenTop}}, for each $\alpha \in J$. 

We turn to show the universality property. Let us have 
for some $(Z, \lambda )$  morphisms
$\langle g_{\alpha } : (Y_{\alpha }, \nu _{\alpha }) \to (Z, \varrho )
\mid \alpha \in J \rangle $ (a sink in {\bf{GenTop}}) 
such that 
$$
Ug_{\alpha }=h \psi _{\alpha } {\text{ for all }} \alpha \in J
$$ 
for some $h: X \to UZ$. Then
$$
\nu _{\alpha } \supset g_{\alpha } ^{-1} ( \varrho ) = (U g_{\alpha }) 
^{-1} ( \varrho ) =( h \psi _{\alpha } ) ^{-1} ( \varrho ) 
= \psi _{\alpha } ^{-1}
h^{-1} ( \varrho )\,.
$$
Hence
$$
\forall \alpha \in J \,\,\,\,
h^{-1} ( \varrho ) \subset \{ M \subset X \mid \psi _{\alpha }^{-1} (M) \in \nu
_{\alpha } \} 
$$
i.e.,
$$
h^{-1} ( \varrho ) \subset \{ M \subset X \mid \forall \alpha \in J \,\,\,\,
\psi _{\alpha }^{-1} (M) \in \nu _{\alpha } \} \,. 
$$
Thus $h = U h'$ for a continuous map $h': (X, \mu ) \to (Z, \varrho ) $.
\end{proof}


We recall \ref{3.1} for notations in the following proofs.


\vskip.2cm

{\bf{Proposition 3.4.}} {\it{Let $X$ be a set. Let us have a source 
$\langle \varphi _{\alpha } :
X \to U(Y_{\alpha }, d_{\alpha })= Y_{\alpha } \mid \alpha \in J \rangle $ 
in {\bf{Set}}. 
Then its initial lift (i.e., weak structure) in {\bf{GenTop}} 
is $(X, c)$, where for $A \subset X$ we
have 
$$
cA:=\cap _{\alpha \in J} 
\varphi _{\alpha } ^{-1} d_{\alpha }\varphi _{\alpha } (A).
$$}}


\begin{proof}
Obviously $c$ is increasing and monotonous. To show idempotence
of $c$, it is sufficient to show $c^2(A) \subset c(A)$ for each $A \subset
X$, i.e.,
$$
\cap _{\beta \in J} 
\varphi _{\beta } ^{-1} d_{\beta }\varphi _{\beta }
[\cap _{\alpha \in J}
\varphi _{\alpha } ^{-1} d_{\alpha }\varphi _{\alpha } (A)] 
\subset
\cap _{\beta \in J} 
\varphi _{\beta } ^{-1} d_{\beta }\varphi _{\beta } (A).
$$
It suffices to show the inclusion (\ref{3.4.1}) 
obtained by deleting here $\cap _{\beta \in J}$
from both sides. We have
\begin{equation}\label{3.4.1}
\begin{array}{l}
\varphi _{\beta } ^{-1} d_{\beta }\varphi _{\beta }
[\cap _{\alpha \in J}
\varphi _{\alpha } ^{-1} d_{\alpha }\varphi _{\alpha } (A)] 
\subset
\varphi _{\beta } ^{-1} d_{\beta }\varphi _{\beta }
[\varphi _{\beta } ^{-1} d_{\beta }\varphi _{\beta } (A)] \\
\subset
\varphi _{\beta } ^{-1} d_{\beta }
d_{\beta }\varphi _{\beta } (A)
=
\varphi _{\beta } ^{-1} d_{\beta } \varphi _{\beta } (A) \,.
\end{array}
\end{equation}
Here we used $[\cap _{\alpha \in J}
\varphi _{\alpha } ^{-1} d_{\alpha }\varphi _{\alpha } (A)] \subset 
\varphi _{\beta } ^{-1} d_{\beta }\varphi _{\beta } (A)$ and $\varphi _{\beta }
\varphi _{\beta } ^{-1} (B_{\beta }) \subset B_{\beta }$ for $B_{\beta } 
\subset Y_{\beta }$ and $d_{\beta }^2= d_{\beta }$. Thus we have obtained the
claimed inclusion \ref{3.4.1}.

Let $A \subset X$. Then 
$$
cA= \cap _{\alpha \in J} \varphi _{\alpha } ^{-1} d_{\alpha } \varphi _{\alpha
} (A) 
\subset 
\varphi _{\alpha } ^{-1} d_{\alpha } \varphi _{\alpha } (A)
$$ 
shows that $\varphi _{\alpha } :
X \to Y_{\alpha }$ becomes (underlies) a continuous map
$f_{\alpha }: (X,c) \to (Y, d_{\alpha })$ for each $\alpha \in J$.

We turn to show the universality property. Let us have 
for some $(Z, e )$  morphisms
$\langle g_{\alpha } : (Z, e ) \to (Y_{\alpha }, d _{\alpha }) \mid \alpha \in
J \rangle $ (a source in {\bf{GenTop}}) 
such that 
$$
Ug_{\alpha }=\varphi _{\alpha } h\, {\text{ for all }} \alpha \in J
$$ 
for some $h:Z \to X$. We have to show that $h$ becomes (underlies) a
continuous map $h':(Z,e) \to (X,c)$. That is, we have to show for $C \subset
Z$ that 
$$
eC \subset h^{-1} c h C = h^{-1}  \cap _{\alpha \in J} \varphi _{\alpha }
^{-1} d_{\alpha } \varphi _{\alpha } h C =
h^{-1} ( \cap _{\alpha \in J} \varphi _{\alpha } ^{-1} d_{\alpha }  
\varphi _{\alpha } ) h C \,.
$$

For $C \subset Z$ we have, for each $\alpha \in J$, 
by continuity of $g_{\alpha }$ that 
\[
\begin{array}{c}
eC \subset g_{\alpha } ^{-1} d_{\alpha }g_{\alpha } C = 
(Ug_{\alpha }) ^{-1} d_{\alpha } (Ug_{\alpha }) C = \\
(\varphi _{\alpha } h) ^{-1} d_{\alpha } (\varphi _{\alpha } h) C
= h^{-1} (\varphi _{\alpha } ^{-1} d_{\alpha } \varphi _{\alpha }) h C\,.
\end{array}
\]
Therefore we have
\[
\begin{array}{c}
eC \subset \cap _{\alpha \in J} 
h^{-1} (\varphi _{\alpha } ^{-1} d_{\alpha } \varphi _{\alpha }) h C =
h^{-1} ( \cap _{\alpha \in J}
\varphi _{\alpha } ^{-1} d_{\alpha } \varphi _{\alpha } ) h C =
h^{-1} c h C\,,
\end{array}
\]
which shows that $h = U h'$ for a continuous map $h':(Z,e) \to (X,c)$.
This shows the universality property of $(X,c)$ and ends the proof of the 
proposition.
\end{proof}


For the next proposition we will need the following concept.

\vskip.2cm

{\bf{Definition 3.5.}} (For the special case of {\bf{PrTop}} \cite{CFK},
Ch. 16B, Topological modifications, in full generality
\cite{DT}, pp. xiv-xv, Ch. 4.6, Idempotent hull and weakly hereditary core,
and \cite{Ca}, p. 74, Definition 6.9 and p. 87, Proposition 7.6.)
Let $(X, \gamma )$ be a CS (with $\gamma $ not necessarily idempotent). 
Let $\lambda $ be any ordinal (also $\kappa $ will denote here ordinals)
and $A \subset X$.
Then $\gamma ^{\lambda }(A)$ is defined in the following way. We let
$
\gamma ^0(A):=A.
$ 
For $\lambda =\kappa +1$ we let 
$
\gamma ^{\lambda }(A):=\gamma \left( \gamma ^{\kappa }(A) \right) .
$
For $\lambda $ a limit ordinal we let 
$
\gamma ^{\lambda }(A):= \cup _{\kappa < \lambda} \gamma ^{\kappa }(A) .
$
Then for any $A \subset X$ there is a smallest ordinal $\lambda _0$ such that
$
\gamma ^{\lambda _0+1} (A) = \gamma ^{\lambda _0} (A) .
$
(Clearly $\lambda _0$ has a cardinality at most $|X|$.)
Then the operator $A \mapsto \gamma ^{\infty } (A):=\gamma ^{\lambda _0} (A)$ 
is called the
{\it{idempotent hull (}}or {\it{transfinite iteration)
of}}\, $\gamma $. Clearly $\gamma ^{\infty }$ is
increasing, monotonous and idempotent. Namely, it equals $\gamma ^{\lambda
_1}$, where $\lambda _1$ is the initial ordinal of the cardinal successor of
$|X|$ for $X$ infinite, or $\lambda _1 = \omega $ for $X$ finite,
for which these properties are obvious.

\vskip.2cm

We will write id for the {\it{identity operation}} (its domain will be clear
from the context).

\vskip.2cm


{\bf{Proposition 3.6.}} {\it{Let $X$ be a set. Let us have a sink $\langle \psi
_{\alpha }: U(Y_{\alpha }, d _{\alpha })= Y_{\alpha } \to X \mid \alpha \in J
\rangle $ 
in {\bf{Set}}. Then its final lift (strong structure) in {\bf{GenTop}}  
is $(X, c )$, where 
\[
\begin{array}{c}
c {\text{ is the idempotent hull }} 
\gamma ^{\infty } 
{\text{ of the}} 
{\text{ closure operator $ \gamma $}} \\
{\text{on }} P(X) 
{\text{ given by}} 
\gamma A := \left[ \cup _{\alpha \in J} 
\psi _{\alpha } d_{\alpha } \psi _{\alpha } ^{-1} (A) \right] \cup A.
\end{array}
\]
}}


\begin{proof}
Clearly $ \gamma $ is
increasing and monotonous. Hence its idempotent hull is increasing, 
monotonous and idempotent. Let $A \subset X$. Then 
\[
\begin{array}{c}
cA = \left[ [\cup _{\alpha \in J} 
\psi _{\alpha } d_{\alpha } \psi _{\alpha } ^{-1} ] \cup  {\text{{\rm{id}}}}
\right] ^{\infty } (A)
\supset [\cup _{\alpha \in J} 
\psi _{\alpha } d_{\alpha } \psi _{\alpha } ^{-1} (A)] \cup A \\
\supset
\cup _{\alpha \in J} 
\psi _{\alpha } d_{\alpha } \psi _{\alpha } ^{-1} (A) \supset
\psi _{\alpha } d_{\alpha } \psi _{\alpha } ^{-1} (A)
\end{array}
\]
shows that $\psi _{\alpha }$ becomes (underlies) a continuous map $f_{\alpha
}: (Y_{\alpha }, d_{\alpha } ) \to (X,c)$, for each $\alpha \in J$.

We turn to the universality property. Let us have 
for some $(Z, e )$  morphisms
$\langle g_{\alpha } : (Y_{\alpha }, d _{\alpha }) \to (Z, e ) \mid \alpha \in
J \rangle $ (a sink in {\bf{GenTop}}) such that 
$$
Ug_{\alpha }=h \psi _{\alpha } {\text{ for all }} \alpha \in J
$$ 
for some $h:X \to Z$. We have to show that $h$ becomes (underlies) a
continuous map $h': (X,c) \to (Z,e) $. 
That is, we have to show for $C \subset
Z$ that 
$$
eC \supset h c h^{-1} C = h  
\left[ [\cup _{\alpha \in J} 
\psi _{\alpha } d_{\alpha } \psi _{\alpha } ^{-1} (A)] \cup {\text{id}} \right]
^{\infty } h^{-1} C \,.
$$

For $C \subset Z$ we have, for each $\alpha \in J$, 
by continuity of $g_{\alpha }$ that 
\[
\begin{array}{c}
eC \supset g_{\alpha } d_{\alpha }g_{\alpha }  ^{-1} C = 
(Ug_{\alpha }) d_{\alpha } (Ug_{\alpha }) ^{-1} C = \\
(h \psi _{\alpha }) d_{\alpha } (h \psi _{\alpha } )  ^{-1} C =
h (\psi _{\alpha } d_{\alpha } \psi _{\alpha } ^{-1} ) h^{-1} C\,.
\end{array}
\]

Therefore we have
\[
\begin{array}{c}
eC \supset \cup _{\alpha \in J} 
h (\psi _{\alpha } d_{\alpha } \psi _{\alpha } ^{-1} ) h  ^{-1} C =
h ( \cup _{\alpha \in J}
\psi _{\alpha } d_{\alpha } \psi _{\alpha } ^{-1}) h ^{-1} C \,,
\end{array}
\]
which implies together with
 $C \subset eC$ that for $\gamma = [\cup _{\alpha \in J}
\psi _{\alpha } d_{\alpha } \psi _{\alpha } ^{-1}] \cup {\text{id}}$  we have
\begin{equation}\label{*}
\begin{array}{l}
h \gamma h^{-1} C=
h \left( [ \cup _{\alpha \in J}
\psi _{\alpha } d_{\alpha } \psi _{\alpha } ^{-1} ] \cup {\text{id}} \right)
h ^{-1} C \\
= h ( \cup _{\alpha \in J}
\psi _{\alpha } d_{\alpha } \psi _{\alpha }  ^{-1} h ^{-1} C) \cup h h^{-1}C \\
\subset \left[ h ( \cup _{\alpha \in J}
\psi _{\alpha } d_{\alpha } \psi _{\alpha }  ^{-1} ) h ^{-1} C \right] \cup C 
\subset eC\,.
\end{array}
\end{equation}
Applying (\ref{*}) to $eC$ rather than $C$, we obtain 
\begin{equation}\label{**}
\begin{array}{l}
h \gamma h^{-1} eC \subset e(eC)=eC \,.
\end{array}
\end{equation}

Now we show by transfinite induction that for each ordinal $\lambda $ and any
$C \subset Z$ we have 
\begin{equation}\label{***}
\begin{array}{l}
h \gamma ^{\lambda } h^{-1} C \subset eC\,.
\end{array}
\end{equation}

For $\lambda =0$ (\ref{***}) is evident 
(since $h h^{-1} C \subset C \subset eC$), and
for $\lambda =1$ this is (\ref{*}). Now let for $\lambda = \kappa +1$ 
(\ref{***}) hold for the ordinal $\kappa $ and
any $C \subset Z$, and
we prove it for $\lambda $ and any $C \subset Z$. We have
\[
\begin{array}{c}
h \gamma ^{\kappa +1} h^{-1} C= 
h \gamma \, {\text{id}} \, \gamma ^{\kappa } h^{-1} C \subset  \\
h \gamma h^{-1} (h \gamma ^{\kappa } h^{-1} C) \subset
h \gamma h^{-1} (eC) \subset eC \,.
\end{array}
\]
In the first inclusion we used that $h^{-1} h$ was increasing, in the second
inclusion we used the induction hypothesis 
and in the third inclusion we used (\ref{**}).  
Now let $\lambda $ be a limit ordinal, and let us suppose that we know
(\ref{***}) for all ordinals 
$\kappa < \lambda $ and for all $C \subset Z$. Then
\[
\begin{array}{c}
h \gamma ^{\lambda } h^{-1} C =
h \cup _{\kappa < \lambda } \gamma ^{\kappa } h^{-1} C =
\cup _{\kappa < \lambda } (h \gamma ^{\kappa } h^{-1} C) \subset
\cup _{\kappa < \lambda } eC = eC \,.
\end{array}
\]
Therefore, by transfinite induction, for each ordinal $\lambda $ we have
(\ref{***}).
In particular, for the ordinal
$\lambda _0$ associated to the set $C$ (cf. Definition 3.5)
and to the operation $\gamma $ we obtain
\[
\begin{array}{c}
h c h ^{-1} C = 
h \gamma ^{\infty } h ^{-1} C = h \gamma ^{ \lambda _0} h ^{-1} C
\subset eC\,, 
\end{array}
\]
which shows that $h = U h'$ for a continuous map $h':(X,c) \to (Z,e)$.
This shows the universality property of $(X,c)$, and ends the proof of the 
proposition.
\end{proof}


\begin{proof} The proof of 
Theorem 3.1 follows from any of Propositions 3.2, 3.3, 3.4, 3.6.
\end{proof}


{\bf{Remark 3.7.}} In Proposition 3.6 it 
is in fact necessary to use transfinite iteration over all ordinals.
As an example, let us take the set $\{ 0,1 \} $, and for some
set $J$ we consider $X := \{ 0,1 \} ^J $. We consider the following sink
to $X$. We let 
\[
\begin{array}{c}
\{ Y_{\alpha } \mid \alpha \in J \}= \{ Y \subset X \mid
|Y|=2, \,\,Y= \{ y_1,y_2 \}, \\ 
{\text{and }} y_1 {\text{ and }} y_2 
{\text{ differ in exactly one coordinate}} \} \,.
\end{array}
\] 
On each $Y_{\alpha }$ we consider the indiscrete
topology, i.e., $d _{\alpha }\emptyset = \emptyset $ and the closure of a
non-empty subset is $Y_{\alpha }$. 
The maps $\psi _{\alpha }$ are the natural injections $Y
_{\alpha } \to X$. We investigate the strong structure on $X$ associated to the
sink $\langle \psi _{\alpha }: Y_{\alpha } \to X \mid \alpha \in J \rangle $. 
Let $x_0 \in X$ be the point with all
coordinates $0$. We write for $A \subset X$
$$
\gamma A := [ \cup _{\alpha \in J} \psi _{\alpha } d _{\alpha } \psi _{\alpha }
^{-1} A] \cup A\,.
$$
Then
$
\gamma \{ x_0 \} = 
\{ x \in X \mid x {\text{ has at most one non-zero coordinate}} \} \,.
$   
Similarly,
$
\gamma ^2 \{ x_0 \} = 
\{ x \in X \mid x {\text{ has at most two non-zero coordinates}} \} \,, 
$
and, in general, for any ordinal $\lambda $,
$
\gamma ^{\lambda } \{ x_0 \} = 
\{ x \in X \mid x {\text{ has at most }} | \lambda | {\text{ non-zero 
coordinates}} \} \,. 
$ 
Therefore the smallest ordinal $\lambda _0$ such that $\gamma ^{ {\lambda }_0
+1 } \{ x_0 \} = \gamma ^{\lambda _0} \{ x_0 \} $ is the 
initial ordinal belonging to the cardinal $|J|$, that can be arbitrarily
large.

\vskip.2cm


Now we give some corollaries to Theorem 3.1 and Propositions 3.2, 3.3, 3.4,
3.6, to some particular kinds of 
initial and final lifts (weak and strong structures), 
and of limits and colimits.

If the source consists of a single map $\varphi : X \to U(Y, \nu) =Y$, 
and $\varphi $ is injective, then the weak structure is called
{\it{subspace of $(Y, \nu )$}}. Here we may suppose that $\varphi $ 
{\it{is actually
the embedding of a subset}}, which we will suppose. If the sink 
consists of a single map $\psi : U(Y, \nu) =Y \to X$, 
and $\psi $ is surjective, then the strong structure is called
{\it{quotient of}} $(Y, \nu )$. 

For $(Y, \nu )$ a GTS and $X \subset Y$, we write $(Y, \nu ) | X := (X, \nu
|X)$, where for ${\cal{A}} \subset P(X)$ we have 
${\cal{A}}|X=\{ A \cap X \mid A \in {\cal{A}} \} $, \S 1, {\bf{2}} 
(cf. also \cite{Sa13}).


\vskip.2cm

{\bf{Corollary 3.8.}} (For supratopologies cf. \cite{KM}, Proposition 4.10, or 
\cite{KM}, Theorem 4.4 and Proposition 4.9. Cf. also \cite{Sa13}.) 
{\it{Subspaces exist in the category {\bf{GenTop}}, and they
can be given as follows. For $(Y, \nu )$ and $i:X \to Y$ an embedding of a
subset, the subspace structure on $X$ is $(Y, \nu )|X = (X, \nu |X)$.
For $(Y,d)$ and $i:X \to Y$ an embedding of a
subset, the subspace structure $(X,c)$ on $X$ is given by
$c(A) :=d(A) \cap X$ for $A \subset X$.
$ \square $}}
 
\vskip.2cm


{\bf{Corollary 3.9.}} 
{\it{Quotient spaces exist in the category {\bf{GenTop}}, 
and
they can be given as follows. For $(Y, \nu )$ and an onto map $q:Y \to X$
in the category \,{\bf{Set}}
the quotient space structure on $X$ (by the map $q$)
is $(X, \{ M \subset X \mid q^{-1}(M) \in \nu \} $. 
For $(Y,d)$ and an onto map $q:Y \to X$ 
in the category \,{\bf{Set}}
the quotient space structure on $X$ (by the map $q$) is $(X,c)$, where 
$c$ is the idempotent hull
of the closure operator on $P(X)$ given by $A \mapsto  
qdq ^{-1}(A) \cup A$. 
$ \square $}}

\vskip.2cm


{\bf{Corollary 3.10.}} (For supratopologies cf. \cite{KM}, Proposition 4.10,
or \cite{KM}, Theorem 4.4 and Proposition 4.9.)
{\it{Products exist in the category {\bf{GenTop}}, and they
can be given as follows. For $\langle (Y_{\alpha }, \nu _{\alpha }) \mid
\alpha \in J \rangle $ 
the product is
$(\prod _{\alpha \in J} Y_{\alpha }, \mu )$ 
(together with the natural projections 
$\pi_{\alpha }$), where $\mu := \{ \cup _{\alpha \in J} \pi _{\alpha }^{-1}
(M_{\alpha }) \mid M_{\alpha } \in \nu _{\alpha } \} $.
For $\langle (Y_{\alpha }, d _{\alpha }) \rangle \mid \alpha \in J \rangle $ 
the product is $(\prod _{\alpha \in J} Y_{\alpha }, c )$,
where for $M \subset \prod _{\alpha \in J} X_{\alpha }$ we have $c(M) := 
\prod _{\alpha \in J} d_{\alpha } \pi _{\alpha} M$. 
$ \square $}}


\vskip.2cm

{\bf{Definition 3.11.}}
The {\it{Cs\'asz\'ar 
product of GTS's}} is defined as follows. Let $J$ be an
index set, let $Y_\alpha $ for $\alpha\in J$ be sets, and $X=\prod _{\alpha \in
J}Y_{\alpha }$. Suppose that, for $\alpha \in J$, $\nu_ {\alpha }$ is a GT on
$Y_\alpha$. 
Let ${\cal{B}} := \{ \prod _{\alpha \in J} N_{\alpha } \mid 
N_{\alpha } \in \nu _{\alpha } $ and, with the exception
of finitely many indices $\alpha $, 
$N_{\alpha }=M_{\nu _{\alpha }} \} $, where
$M_{\nu _{\alpha}} := \cup_{B \in \nu _{\alpha }} B$. 
The GT on $X$ having ${\cal{B}}$ as a base 
is called the Cs\'asz\'ar product of the GTS's $\langle (Y_{\alpha }, \nu
_{\alpha }) \mid \alpha \in J \rangle $.

\vskip.2cm


{\bf{Remark 3.12.}} Obviously, Cs\'asz\'ar's products of strong GTS's 
are finer than the product GTS's (they have the same underlying set $\prod
_{\alpha\in J} Y_{\alpha }$), but if strongness is omitted, they are
in general incomparable, even for $|J|=2$. Also
the (categorical) product of GT's in general is not a topology, even if the
factors are topological spaces, but the Cs\'asz\'ar product is the topological
product in this last case.

Categorical products coincide with the Cs\'asz\'ar product 
only in some particular cases. The empty Cs\'asz\'ar product is $(X_0, \{
\emptyset , X_0 \} )$ with $|X_0|=1$, while the empty product is $(X_0, \{ 
\emptyset \} )$ with $|X_0|=1$. For $|J|=1$ both the Cs\'asz\'ar product and
the product equal the unique factor. If some $Y_{\alpha }$ is empty, both the
Cs\'asz\'ar product and the product are $( \emptyset , \{ \emptyset \} )$.
For $|J| \ge 2$ and $\forall \alpha \in J\,\,\,\,Y_{\alpha } \ne \emptyset $,
the
equality of the Cs\'asz\'ar product and the product is equivalent to that each
$(Y_{\alpha }, \nu _{\alpha })$ is strong, and there is at most one $\alpha \in
J$ such that $\nu _{\alpha } \ne \{ \emptyset , Y_{\alpha } \} $.

\vskip.2cm


{\bf{Corollary 3.13.}} {\it{Sums exist in the category {\bf{GenTop}}, and they
can be given as follows. For $\langle (Y_{\alpha }, \nu _{\alpha }) \mid
\alpha \in J \rangle $ 
the sum is
$(\coprod _{\alpha \in J} Y_{\alpha }, \mu )$ 
(together with the natural injections $\iota _{\alpha }$),
where $\mu := \{
\cup _{\alpha \in J }  \iota _{\alpha } N_{\alpha} \mid \forall \alpha \in J\,\,
\,\,N_{\alpha } \in \nu _{\alpha } \} $.
For $\langle (Y_{\alpha }, d _{\alpha }) \mid \alpha \in J \rangle $ 
the sum is $(\coprod _{\alpha \in J} Y_{\alpha }, c )$,
where for 
$N_{\alpha } \subset Y _{\alpha }$
we have $c(\cup _{\alpha \in J} \iota _{\alpha } N_{\alpha }):=
\cup _{\alpha \in J} \iota _{\alpha }
d_{\alpha }N_{\alpha }$.
$ \square $}}


\vskip.2cm

{\bf{Corollary 3.14.}} (For supratopologies cf. \cite{KM}, Proposition 4.1 and
Theorem 4.4.)
{\it{Let $X$ be a set. Then all generalized topologies on
$X$ form a complete lattice, with $(X, \mu ) \le (X, \nu )$ meaning that
the identical map of $X$ is (underlies) a continuous map $ (X, \nu ) \to    
(X, \mu ) $. The union, or intersection 
of a set of generalized topologies 
$\{ (X, \nu _{\alpha } ) \mid \alpha \in J \} $ is $(X, \mu )$, where
$\mu := \{ \cup _{\alpha \in J} N_{\alpha } 
\mid N_{\alpha } \in \nu _{\alpha } \}
$, or $\mu := \cap _{\alpha \in J}\nu _{\alpha }$, respectively.
The union, or intersection 
of a set of generalized topologies $\{ (X, d _{\alpha }) \mid {\alpha } \in J
\} $ 
is $(X, c )$, where
$c(A) := \cap _{\alpha \in J} d_{\alpha } (A)$ for $A \subset X$, or
$c$ is the idempotent hull of the closure operator
$A \mapsto [\cup _{\alpha \in J} d_{\alpha } (A)] \cup A$ for $A \subset X$,
respectively.
$ \square $}}

\vskip.2cm


In the last formula ``union with $A$'' is necessary only for $J = \emptyset $.


\section{$T_{3.5}$, normal, compact, Lindel\"of GTS's and the analogue of
Tychonoff's embedding theorem}


Further on, products of GTS's will be meant in the categorical sense, 
cf. \S 2, {\bf{4}}. Cs\'asz\'ar's products will not be used.

Cs\'asz\'ar \cite{Cs2} introduced a useful {\it{GT on the set of real numbers
$\mathbb{R}$}} as follows. It has as a base 
$$
\beta := \{ (-\infty , s) \mid s \in {\mathbb{R}} \} \cup
\{ (t, \infty ) \mid t \in {\mathbb{R}} \} \,.
$$
This is a strong GT.

Henceforth, {\it{we assign the notation $\gamma$
just for this GT}}. We believe that this GTS is the appropriate choice for
$\mathbb{R}$ as a GTS. Indeed, $({\mathbb{R}}, \gamma )$ as a GTS has a similar
role as the 
standard topology on $\mathbb{R}$ in general topology. 
Similarly, {\it{we use the notation $([0,1],\gamma _0)$ for the subspace
$[0,1]$ of
$\mathbb{R}$}} (i.e., $\gamma _0 = \gamma | [0,1]$, cf. Corollary 3.8).

\vskip.2cm

{\bf{Remark 4.1.}} (\cite{A}, Remark after Example 2.4.)
The GTS $(\mathbb{R}, \gamma )$ is $T_4$. (Namely a simple
discussion shows that its closed sets are just the 
convex subsets of $\mathbb{R}$ which are closed in the usual sense.
Then two non-empty disjoint closed sets can
be included into two disjoint halflines which are open in the usual sense. 
The $T_1$ property is evident.) 
By \cite{Cs2}, Proposition 2.5 normality is closed-hereditary, as well as the
$T_1$ property, hence 
$T_4$ property of $({\mathbb{R}}, \gamma )$ implies the $T_4$  property of
$([0, 1],\gamma _0)$.

\vskip.2cm

{\bf{Remark 4.2.}} Let $n \ge 2$. We suppose that 
the appropriate GTS analogue of the $n$-dimensional real vector space is not
the $n$'th power of $(\mathbb{R}, \gamma )$ (this depends on some 
preassigned representation of the $n$-dimensional real vector space 
as a direct sum of $1$-dimensional subspaces). 
We suppose the proper choice should be the GT with base all
open (in the usual sense) halfspaces. Then the closed sets are exactly the
closed (in the usual sense) convex subsets, and the associated closure operator
is the closed (in the usual sense)
convex hull of a subset (cf. e.g., \cite{BF}, \S 1, 3, 
\cite{DS}, V.2.7 Theorem 10,
\cite{Sch}, Theorem 1.3.7). These have been 
extensively investigated in
geometry and functional analysis, cf. e.g. the just cited three books.
(For weak topologies on locally convex topological vector spaces there is an
analogue of this construction: we consider inverse images of $(r, \infty )$,
where $r \in {\mathbb{R}}$, by
continuous linear functionals, as basic generalized
open sets, and then the generalized closed sets are
still the closed (in the usual sense) convex sets, and the associated closure
operator is the closed convex hull of a subset.
For this, including the
necessary definitions, cf. \cite{DS}.)

Observe that this space is not normal for $n \ge 2$: the closed sets $F_1:=
x_1 \dots
x_{n-1}$-coordinate hyperplane and $F_2 := \{ (x_1, \dots , x_n) \mid \forall
i \in \{ 1, \ldots , n \} \,\,x_i >
0,\,\,x_1 \dots x_n \ge 1 \} $ 
are disjoint closed subsets which cannot be included into
two open subsets. In fact, non-empty disjoint open sets $M_1,M_2$ 
are unions of open (in the usual sense) half-spaces $ \{ M_{1, \alpha } \mid
\alpha \in A \} $ and
$\{ M_{2, \beta } \mid \beta \in B \} $. 
Then each $M_{1, \alpha }$ and each $M_{2, \beta }$,
being disjoint, have
parallel boundary hyperplanes. So $M_1,M_2$ also are open (in the usual sense)
half-spaces, having parallel boundary hyperplanes. Therefore
we may suppose that $|A|=|B|=1$, $A =  \{
\alpha _0 \} $, $B = \{ \beta _0 \} $, and the boundary hyperplanes of $M_{1,
  \alpha _0}$ and $M_{2, \beta _0}$ (in the usual sense) are parallel to $F_1$.
If $F_1
\subset M_{1, \alpha _0}$, then $M_{1, \alpha _0}$ contains a parallel slab 
containing some usual open $\varepsilon $-neighbourhood $F_{1, \varepsilon }$ 
of $F_1$, and then $M_1 \cap M_2 \supset F_{1, \varepsilon } \cap F_2 \ne
\emptyset $. 

\vskip.2cm

By the GTS $([0,1], \gamma _0)$ 
\cite{Cs2} exhibited a generalization of the Urysohn lemma for
 normal GTS's, as follows. 


\vskip.2cm

{\bf{Theorem 4.3}}
(\cite{Cs2}, Theorem 3.3)
 {\it{Let 
$(X, \mu )$ be a normal GTS, and $F, F' \subset X$ be disjoint $\mu$-closed
 sets. Then there exists a continuous function $f:(X, \mu )
 \rightarrow ([0,1], \gamma _0)$ such that $f(x)=0$ for $x\in F$ and 
$f(x)=1$ for $x\in F'$. $\square $}} 

\vskip.2cm


By the above discussions we are able to define and investigate 
$T_{3.5}$ spaces.


\vskip.2cm

{\bf{Definition 4.4.}}
A GTS $(X, \mu )$ is {\it{completely regular}} if 
for each $x \in X$ and each $\mu$-closed set $F$ of $X$ not
containing $x$, there is a continuous function $f:(X, \mu )
\rightarrow ([0, 1], \gamma _0)$ 
such that $f(x)=0$ and $f(F) \subset \{ 1 \} $. We say that
$(X, \mu )$ is $T_{3.5}$ (or
{\it{Tychonoff}}) if it is a completely regular $T_1$ space.
 
\vskip.2cm

We will write $T_{3.5}$ spaces rather than Tychonoff spaces.
In Definition 4.4 we can write $f(F)=\{ 1 \} $ if we require $F \ne \emptyset $,
which we may
suppose. (For $F = \emptyset $ $(X, \mu )$ is strong, and then we can take the
identically $0$ function to $([0, 1], \gamma _0)$, which is now continuous.)
In particular, $T_{3.5}$ is equivalent to completely regular $T_0$,
since by \S 2, {\bf{6}}, a regular $T_0$ space is $T_1$.
Clearly, every $T_{3.5}$, or completely regular space is a $T_3$, or regular
space (cf. \S 2, {\bf{6}}), moreover, by Theorem 4.3, every
$T_4$ space is a $T_{3.5}$ space.
Also, as in the case of regularity (cf. \S 2, {\bf{6}}), 
a GTS of the form $(X, \{
\emptyset \} ) $ is vacuously completely regular, but complete regularity of a
GTS $(X, \mu )$ with $\mu \ne \{ \emptyset \} $ 
implies its regularity, and hence its strongness (cf. \S 2, {\bf{6}}).


We are
going to exhibit a generalized version of Tychonoff's embedding theorem to
obtain a necessary and sufficient condition for a GTS to be a
$T_{3.5}$ space.


\vskip.2cm

{\bf{Example 4.5.}}
Consider the GTS $([0, 1], \gamma _0 )$. Then
\[
\begin{array}{c}
\gamma _0 = \{ \emptyset, \left[ 0, 1 \right] \} \cup 
\{ \left[ 0, p) \mid p \in (0,1 \right] \} \cup \\
\left. \{ (q,1 \right] \mid q \in \left[ 0,1) \} \cup 
\{ \left[ 0, r) \cup (s, 1 \right] \mid r, s \in (0, 1),\,\,r \le s \}
\,. \right.
\end{array}
\]
By Remark 4.1 $([0, 1],\gamma _0)$ is $T_4$, hence it 
is a $T_{3.5}$ GTS as well. 

\vskip.2cm


{\bf{Definition 4.6.}} A source $\langle 
f_{\alpha }: X \to Y_{\alpha } \mid \alpha \in J \rangle $ 
in {\bf{GenTop}}, {\bf{Set}} 
is {\it{point-separating}} (also called {\it{monosource}}) if
for $x_1,x_2 \in X$ distinct there is an $\alpha \in J$ such that
$f_{\alpha }(x_1) \ne f_{\alpha }(x_2)$.


\vskip.2cm

{\bf{Proposition 4.7.}} (For supratopologies, for regularity cf. Theorem 4.4.)
{\it{Let $\langle f_{\alpha } : X \to U(Y_{\alpha }, \mu
_{\alpha }) = Y_{\alpha } \mid \alpha \in J \rangle $ 
be a source in {\bf{Set}}.\\ 
(1)\ If each $(Y_{\alpha }, \mu_{\alpha })$ is regular or
completely regular,
then also the initial (weak) GTS structure on $X$ defined by this source is 
regular or completely regular.\\
(2) If each $(Y, \mu_{\alpha })$ is $T_0$, 
$T_1$, $T_2$, $T_3$ or $T_{3.5}$, and 
still $\langle f_{\alpha } \mid \alpha \in J \rangle $ is point-separating,
then also the initial (weak) GTS structure on $X$ defined by this source is 
$T_0$, $T_1$, $T_2$, $T_3$ or $T_{3.5}$.
In particular, subspaces and products of $T_0$,
$T_1$, $T_2$, regular, $T_3$, completely regular or $T_{3.5}$ GTS's are $T_0$,
$T_1$, $T_2$, regular, $T_3$, completely regular or $T_{3.5}$.}}


\begin{proof}
In case (1) we give the proof for complete regularity, and in case (2) for
$T_1$. All other proofs are analogous.
 
We prove (1) for complete regularity. 
Suppose we have a source $\langle \varphi _{\alpha }: X \to U(Y_{\alpha }, \nu
_{\alpha })=Y \mid \alpha \in J \rangle $ in {\bf{Set}}, 
with all $ (Y_{\alpha }, \nu _{\alpha })$ completely regular.
By Proposition 3.2
the weak structure w.r.t. this source is the GT on $X$ having as base $\cup
_{\alpha \in J} \varphi _{\alpha }^{-1} (\nu _{\alpha })$. Therefore it
suffices to prove that for $x \in X$, $x \in N:=\varphi _{\alpha }^{-1} 
(N _{\alpha })$, $N _{\alpha } \in \nu _{\alpha }$, $ \alpha \in J$
and $F=X \setminus N$
there is a continuous function $h: (X, \mu ) \to ([0,1], \gamma _0)$ 
such that $h(x)=0$ and $h(F) \subset \{ 1 \} $, i.e., $h^{-1}[0,1) \subset N$. 

Observe that $\varphi _{\alpha } (x)   \in N _{\alpha }$.
We have that $\varphi _{\alpha }$
becomes (underlies) a continuous map $f_{\alpha } : (X, \mu ) \to (Y_{\alpha
}, \nu _{\alpha })$ in {\bf{GenTop}}.
By complete regularity of $(Y_{\alpha }, \nu _{\alpha })$ there
is a continuous function $g_{\alpha }:(Y_{\alpha }, \nu _{\alpha }) 
\to ([0,1], \gamma _0)$ such that 
$g_{\alpha } \left( f _{\alpha }(x) \right) =0$ 
and $g_{\alpha }^{-1}[0,1) \subset N_{\alpha }$. Then we define $h:=g_{\alpha
} \circ f _{\alpha }$ which satisfies the claimed properties.

We prove (2) for $T_1$. Suppose that we have a source like above,
which is additionally point-separating. We use the notation $f_{\alpha }$ like
above. Let $x_1 \ne
x_2$ be in $X$. Then there exists an index $\alpha \in J$ such that $f
_{\alpha } (x_1) \ne f _{\alpha } (x_2)$. Then by the $T_1$ property of
$Y_{\alpha }$ there is some open set $N_{\alpha } $ in $(Y_{\alpha }, \nu
_{\alpha })$ such that $f _{\alpha } (x_1) \in N_{\alpha } \not\ni 
f _{\alpha } (x_2)$. Then $x_1 \in \varphi _{\alpha } ^{-1} N_{\alpha }
\not\ni x_2$, and by continuity of $f_{\alpha }$ we have
$f_{\alpha } ^{-1} N_{\alpha } \in \mu $, proving
the claimed $T_1$ property of $(X, \mu )$.
\end{proof}


{\bf{Corollary 4.8.}} {\it{The GTS on the $n$-dimensional real vector space 
in Remark 4.2 is $T_{3.5}$.}}


\begin{proof}
This GTS is the weak structure w.r.t. all non-$0$ linear
functionals, as functions to the GTS $({\mathbb{R}}, \gamma )$. 
\end{proof}


Now, we are ready to prove the following variant of Tychonoff's embedding
theorem that characterizes $T_{3.5}$ GTS's. As for topological spaces,  
{\it{a map}} $f : (X, \mu ) \to (Y, \nu )$ {\it{in}}\, {\bf{GenTop}}
{\it{is dense}}, if $Z:=f(X)$ {\it{is dense in}} $(Y, \nu )$, i.e., 
if the closure of $Z$ in $(Y, \nu )$ equals $Y$.


We begin with the analogue of the embedding lemma in topology. We call $f:(X,
\mu ) \to (Y, \nu )$ in {\bf{GenTop}} {\it{open}}, if $f( \mu ) \subset \nu $.

\vskip.2cm

{\bf{Lemma 4.9.}} 
{\it{Let $(X,\mu )$ be a GTS, and let us have a point-separating 
source (monosource) $\langle
f _{\alpha } : (X, \mu ) \to (Y_{\alpha }, \nu _{\alpha }) \mid \alpha \in J
\rangle $ in {\bf{GenTop}}. 
Suppose that for any $x \in M \in \mu $ there exists an $\alpha \in J$ and
$f_{\alpha } (x) \in N_{\alpha } \in \nu _{\alpha }$ 
such that $(x \in )\,\,f_{\alpha } ^{-1}(N_{\alpha })
\subset M$. Then the mapping $f:(X, \mu ) \to \prod _{\alpha \in J}
(Y_{\alpha }, \nu _{\alpha
})$ defined by $f(x)= \langle f_{\alpha } \mid \alpha \in J \rangle $ 
satisfies that 
$f$ is a homeomorphism to its image $ \left( \prod _{\alpha \in J}
(Y_{\alpha }, \nu _{\alpha }) \right) | f(X) $.}}


\begin{proof} The proof 
is analogous to the case of topological spaces. Recall that
$\cup _{\alpha
\in J} \pi _{\alpha }^{-1}
(\nu _{\alpha })$ is a base for $\prod _{\alpha \in J}(Y_{\alpha }, \nu
_{\alpha })$ (where the $\pi _{\alpha }$'s are the natural projections). 
Also, by hypothesis, $f$ is injective.

The hypothesis means that $\{ f_{\alpha } ^{-1}(N_{\alpha }) \mid \alpha \in
J,\,\,N_{\alpha } \in \nu _{\alpha } \} $ is a base of $\mu $. We factorize $f:
(X, \mu ) \to \prod _{\alpha \in J} (Y_{\alpha }, \nu _{\alpha })$ across the
image $f(X)$ as $(X, \mu ) {\overset{F} \to } 
\left( \prod _{\alpha \in J} \right.$
\newline
$ \left. (Y_{\alpha }, \nu _{\alpha })
\right) | f(X) \hookrightarrow \prod _{\alpha \in J} (Y_{\alpha }, \nu
_{\alpha })$.
We are going to show that $F$ is open. 

It suffices to show that for all $\alpha \in J$, and for all
$N_{\alpha } \in \nu _{\alpha }$ we have that
$f f_{\alpha }^{-1} (N_{\alpha })$ is open in $ \left( \prod
_{\alpha \in J} (Y_{\alpha }, \nu _{\alpha  } ) \right) \mid f(X) $. 
We have evidently
$f f_{\alpha }^{-1} (N_{\alpha }) = \langle f_{\beta } \mid \beta \in J
\rangle f_{\alpha } ^{-1}(N_{\alpha }) \subset f(X) \cap \pi _{\alpha }^{-1}
(N_{\alpha })$. Next we show the converse inclusion. Let $f(x) \in 
f(X) \cap \pi _{\alpha }^{-1} (N_{\alpha })$. Then $f_{\alpha }(x) = 
\pi _{\alpha } \langle f_{\beta } \mid \beta \in J \rangle (x) = \pi _{\alpha }
f(x) \in N_{\alpha }$, hence $x \in f_{\alpha }^{-1} (N_{\alpha })$, and $f(x)
\in f f_{\alpha }^{-1} (N_{\alpha }) $. This shows the converse
inclusion, hence  $f f_{\alpha }^{-1} (N_{\alpha }) = f(X) \cap \pi _{\alpha
} ^{-1} (N_{\alpha })$ is a basic open set in $ \left( 
\prod _{\alpha \in J} (Y_{\alpha }, \nu _{\alpha }) \right) | f(X) $. 

Thus $F : (X, \mu ) \to
\left( \prod _{\alpha \in J} (Y_{\alpha }, \nu _{\alpha } ) \right)
|f(X) $ is open. Thus $F$
is continuous, open and bijective, hence is a homeomorphism.
\end{proof}


{\bf{Definition 4.10.}}
For $(X, \mu )$ a GTS, we define its $T_0$-{\it{reflection}} as follows. 
We define
$x,y \in X$ equivalent, written $x\equiv y$, by $\forall M \in \mu\,\,(x \in M
\Longleftrightarrow y \in M)$. 
This is an equivalence relation, and we let $Y$ be the
quotient of $X$ w.r.t. this equivalence relation. Hence we have an onto map
$q:X\to Y$, where $q(x)$ is the equivalence class of $x$. The {\it{quotient
space structure on}} $Y$ (cf. Corollary 3.9), {\it{say, $(Y, \nu )$, is the
$T_0$-reflection of}} $(X, \mu )$. Then $(Y, \nu )$ is $T_0$.

\vskip.2cm


{\bf{Proposition 4.11}} 
{\it{(1) Definition 4.10 is correct. Also, $(X, \mu )$ is the initial (weak)
structure associated to the one-element source $Uq:X \to U(X, \mu )$ in
{\bf{Set}}. Moreover, the map $q: (X, \mu ) \to (Y, \nu )$ has the
following universality property. If $(Z, \varrho )$ is a $T_0$ GTS, and $f:
(X, \mu ) \to (Z, \varrho )$ is continuous, then there exists a unique
continuous $h:(Y, \nu ) \to (Z, \varrho )$, such that $hq=f$.

(2) With the above notations, we have the following equivalences: $(X, \mu )$ is
regular (completely regular) $\Longleftrightarrow (Y, \nu )$ is 
regular (completely regular) $\Longleftrightarrow (Y, \nu )$ is
$T_3\,\,(T_{3.5})$.}}


\begin{proof}
(1) For $x_1,x_2,x_3 \in X$, with $x_1 \equiv x_2 \equiv x_3$ we have, for
each $M \in \mu $, that $x_1 \in M \Longrightarrow x_2 \in M \Longrightarrow
x_3 \in M$. The converse implication is proved in the same way, hence $x_1
\equiv x_3$.

The statement about the initial structure is evident.

Let $f:(X, \mu ) \to (Z, \varrho )$, where $(Z, \varrho )$ is $T_0$. Then the
equivalence relation on $(Z, \varrho )$ is the finest one (i.e., 
each equivalence class is a singleton).
Then for $z_1 \ne z_2$ we have that $f^{-1} (z_1), f^{-1}(z_2) \subset X$
contain no $x_1,x_2$ such that $x_1 \equiv x_2$. Therefore the equivalence
relation on $(X, \mu )$ (that corresponds to the partition $\{ q^{-1}(y) \mid
y \in Y \} $) is finer than the equivalence relation corresponding to the
partition $\{ f^{-1} (z) \mid z \in Z \} $. Let $C \in \varrho $. Then $f^{-1}
(C) \in \mu $ is a union of some sets $f^{-1}(z)$ with $z \in Z$, hence is a
union of some equivalence classes associated to $\mu $. Therefore $f:(X, \mu
) \to (Z, \varrho )$ factors as $(X, \mu ) \overset{q}\to (Y, \nu ) 
\overset{h}\to (Z, \varrho )$. Unicity of such an $h$ follows since $q$ is
onto.

(2) The first equivalences follow from the definition of $(Y, \nu )$. The
second equivalences follow since $(Y, \nu )$ is $T_0$, which by regularity
imply $T_1$, and hence $T_3\,\,(T_{3.5})$.
\end{proof}


{\bf{Theorem 4.12}}
{\it{A $GTS$ $(X, \mu )$ is $T_{3.5}$ (completely regular)
if and only if it is homeomorphic
to a subspace $Y$ of a power of the GTS $([0, 1],\gamma _0)$ (has the weak
structure w.r.t. a source $\varphi : X \to U([0, 1],\gamma _0)^J= [0, 1]^J$ 
in {\bf{Set}} --- consisting of a single map ---
for some $J$). For $|X| \ge 2$ (for $\mu \not\subset \{ \emptyset , X \} $)
we may even suppose that $Y$ is dense ($\varphi $ is dense).}}


\begin{proof}
{\bf{1.}} We begin with the $T_{3.5}$ property.

Necessity is a direct consequence of Proposition 4.7 (2), because $([0,
1], \gamma _0)$ is a $T_{3.5}$ space.

We turn to sufficiency. The space
$( \emptyset , \{ \emptyset \} )$ is a subspace of any power of 
$([0,1], \gamma _0)$. For $|X|=1$, $(X, \{ \emptyset \} ) \cong 
([0,1], \gamma _0) ^0$ (cf. \S 2, {\bf{4}}), 
and $(X, \{ \emptyset , X \} )$ is a subspace of
$([0,1], \gamma _0)$. So we may suppose $|X| \ge 2$, that by $T_1$ implies
strongness and $ \mu \not\subset \{ \emptyset , X \} $. 

Put $J= \{ (x, M) \mid x \in M \in \mu \} $. By $\mu \ne \{ \emptyset \} $
we have $|J| \ge 1$. Since $X$ is a
completely regular space therefore for every $ \alpha = (x, M) \in J$ 
there exists a
continuous function $f_{\alpha} : (X, \mu ) \rightarrow ([0, 1], \gamma _0)$ 
such that $f_\alpha(x)=0$ and $f_{\alpha}(X \setminus M) 
\subset \{ 1 \} $. 
Then for $M \ne X$ we have $f_{\alpha}(X \setminus M) = \{ 1
\} $, and for $M = X$ by strongness of $(X, \mu )$ we can
choose for $f_{\alpha}$ the constant $0$ function.
Now we define
$f : (X, \mu ) \rightarrow ([0, 1], \gamma _0)^J$ 
by $f(x) := \langle f_{\alpha}(x)) \mid \alpha \in J \rangle $. 
If $x, y \in X$ and
$x \neq y$, then by the $T_1$ property of $(X, \mu )$,
$ \{ y \} $ is a $\mu $-closed set not containing $x$ --- therefore, there
is $\alpha\in J$ such that $f_{\alpha}(x)=0$ and $f_{\alpha}(y)=1$. Hence,
in presence of the $T_1$ property of $(X, \mu )$,
$ \langle f_{\alpha } \mid \alpha \in J \rangle $ is point-separating.

Applying Lemma 4.9 we obtain that 
the map $f$ is a homeomorphism to its image $(Y, \nu )
:= ([0,1], \gamma _0)^J|f(X) $.

Now we prove the addition about denseness. If $|X| \ge 2$ and we have $T_1$,
then $\mu \not\subset \{ \emptyset , X \} $, and $|J| \ge 1$, and
$(Y, \nu )$ is a subspace of $([0,1], \gamma _0)^J$. Then $(Y, \nu )$ 
is a subspace of $\prod _{\alpha \in J} \pi _{\alpha }Y$.
We distinguish the cases whether $M=X$ or $M \ne X$.
For $M=X$ we have that
$\pi _{\alpha }Y = (\{ 0 \} , \{ \emptyset , \{ 0 \} \} ) $.
All such
factors have product (up to isomorphism)
this same space $(\{ 0 \} , \{ \emptyset , \{ 0 \} \}
)$, and we may omit all such factors.  For $M \ne X$ we have that 
$\pi _{\alpha }Y \supset \{ 0, 1 \} $. Hence
$c_{\gamma _0} \pi _{\alpha }Y = [0,1]$, for all non-omitted factors, which
factors form a set $J'$ $( \subset J$), of
cardinality at least $1$ (by $T_1$ and $|X| \ge 2$). Then a homeomorphic
copy $(Y', \nu ')$ of $(Y, \nu )$ 
is contained in the product $([0,1], \gamma _0)^{J'}$
of the non-omitted factors 
--- namely, $(Y', \nu ') := ([0,1], \gamma _0)^{J'}|
\langle f_{\alpha } \mid \alpha \in J' \rangle (X)$ ---
whose product with $( \{ 0 \} , \{ \emptyset , \{ 0 \} \} )$ is, up to
isomorphism, $(Y, \nu )$.
Since for all $\alpha \in J'$ we have $c_{\gamma _0} \pi _{\alpha }Y = [0,1]$,
therefore by Corollary 3.10 the closure of $Y'$ in $([0,1], \gamma _0 )^{J'}$ 
is $[0,1]^{J'}$, i.e., $Y'$ is dense in $([0,1], \gamma _0 )^{J'}$.

{\bf{2.}}
We turn to the complete regularity.

Spaces of the form $(X, \{ \emptyset \} )$ and (for $X \ne \emptyset $) 
$(X, \{ \emptyset , X \} )$ are completely regular, and have the weak
structures w.r.t. the sources consisting of the single map to 
$([0,1], \gamma _0 )^0 \cong ( \{ 0 \} ,  \{ \emptyset \}  )$, 
and the identically
$0$ map to $([0,1], \gamma _0 )$, respectively. Hence we may suppose $\mu
\not\subset \{ \emptyset , X \} $, when the $T_0$-reflection of $(X, \mu )$
has at least two
points, and when the functions $f_{\alpha }: (X, \mu ) \to ([0,1], \gamma _0
)$ for $M \not\in \{ \emptyset , X \} $ satisfy $f_{\alpha }(x)=0$ and
$f_{\alpha }(X \setminus M ) = \{ 1 \} $. Hence $J \ne \emptyset $, and even
$J' \ne \emptyset $ from part {\bf{1}} of this proof.

Necessity follows, since for a source consisting of a single map $\varphi :
X \to U \left( ([0,1], \gamma _0 )^J \right) = [0,1]^J$ we have that 
$ ( [0,1], \gamma _0 )^J|\varphi (X) $ is $T_{3.5}$
by Proposition 4.7 (2). Then the initial lift (weak structure)
for this source 
$\varphi : X \to U\left( ([0,1], \gamma _0 )^J \right) $ 
is the same as that for the source 
$X {\overset{q}\to} \varphi (X)$ $(\hookrightarrow U \left( 
([0,1], \gamma _0 )^J \right) $,
where $\left( \varphi (X), \gamma _0 ^J|\varphi (X) \right) $ is the
$T_0$-reflection of $(X, \mu )$, and $q$ from Definition 4.10
pointwise coincides with $ \varphi $.
By Proposition 4.11 (2), the $T_{3.5}$ property of $\left( \varphi (X),
\gamma _0^J|\varphi (X) \right) $ implies the complete regularity of
$(X, \mu )$.

For sufficiency, let $(X, \mu )$ be completely regular, and let $(Y, \nu )$ be
its $T_0$-reflection (with an onto map $q:X \to Y$ from Definition 4.10), which
is $T_{3.5}$ by Proposition 4.11 (2), hence is a subspace of some power
$([0,1], \gamma _0) )^J$ by part {\bf{1}} of this proof. Let 
$i:(Y, \nu ) \hookrightarrow 
([0,1], \gamma _0 )^J$ be the inclusion. Then, by Proposition 4.11 (1), $(X,
\mu )$ has the weak structure for the source consisting of the single map
$iq:X \to U \left( ([0,1], \gamma _0 )^J \right) $.
\end{proof}


It is interesting to observe that one can completely describe the continuous
maps of $({\mathbb{R}}, \gamma )$, or of $([0,1], \gamma _0)$, to itself.
We do this in greater generality. Recall that for an ordered set $(X, \le )$
the {\it{order topology}} has as subbase $\{ \{ x \in X \mid x < a \} 
\mid a \in X \} \cup \{ \{ x \in X \mid x > b \} \mid b \in X \} $.

\vskip.2cm

{\bf{Definition 4.13.}} 
{\it{Let $(X, \le )$ be an ordered set. Its order GT has as
base  $\{ \{ x \in X \mid x < a \} 
\mid a \in X \} \cup \{ \{ x \in X \mid x > b \} \mid b \in X \} $.}}
We will write these sets as $( - \infty , a)$ and $(b, \infty )$. 
Analogously we use the notations $( - \infty , a]$ and $[b, \infty )$ for
$\{ x \in X \mid x \le a \} $ and $\{ x \in X \mid x \ge b \} $, respectively.
(Observe that $\pm \infty $ are just symbols and not elements of $X$, even if
$X$ happens to have a minimal or a maximal element. Also, we will use the same
symbols for any other ordered space $Y$ as well, but this will not cause
misunderstanding.)

\vskip.2cm

{\bf{Remark 4.14.}}
Letting $(\tilde{X}, \le )$ be the Dedekind completion of $(X, \le )$, the open
sets of the order GT of $X$
are exactly of the form $ \emptyset $, or $X$ (this only for $|X| \ne 1$),
$\{ x \in X \mid x < {\tilde{x}}_1 \} $, or 
$\{ x \in X \mid x > {\tilde{x}}_2 \} $, or
$\{ x \in X \mid x < {\tilde{x}}_1 {\text{ or }} x > {\tilde{x}}_2 \} $, where
$ {\tilde{x}}_1, {\tilde{x}}_2 \in {\tilde{X}}$ with ${\tilde{x}}_1 \le
{\tilde{x}}_2$, but excluding ${\tilde{x}}_1 = {\tilde{x}}_2 \in {\tilde{X}}
\setminus X$.
Thus, the closed sets are of the form $X$, or $\emptyset $ (this only for $|X|
\ne 1$), or 
$ \{ x \in X \mid x \le {\tilde{x}}_1 \} $
or $\{ x \in X \mid x \ge {\tilde{x}}_2 \} $, or
$\{ x \in X \mid {\tilde{x}}_1 \le x \le {\tilde{x}}_2 \} $, where
$ {\tilde{x}}_1,  {\tilde{x}}_2 \in {\tilde{X}}$ with ${\tilde{x}}_1 \le
{\tilde{x}}_2$, but excluding ${\tilde{x}}_1 = {\tilde{x}}_2 \in {\tilde{X}}
\setminus X$.
I.e., for $|X|\ne 1$, the closed sets in $X$ are exactly the {\it{convex}}
subsets (i.e., which contain with any two
of their points the whole interval between them) which are closed in the (finer)
order topology of $(X, \le )$.

For $X = \emptyset $ (with the void ordering) we have that the order GT is 
$(\emptyset, \{ \emptyset \} )$, that is strong. 
For $|X| = 1$ we have that the order GT is
$(X, \{ \emptyset \} )$, that is not strong. For $|X| \ge 2$ the order GT is
strong (since it is $T_1$, cf. \S 2, {\bf{6}}).

As a generalization of our Remark 4.1, \cite{A} Remark after Example 2.4 claimed
that the order GT associated to any ordered set $(X, \le )$ is $T_4$. (A
detailed simple 
proof can be given for $|X| \le 1$ from the preceding paragraph, while
for $|X| \ge 2$ by using the case distinctions from the second preceding
paragraph.)

\vskip.2cm

{\bf{Proposition 4.15}} 
{\it{Let $(X, \le )$ and $(Y, \le )$ be two ordered sets,
with order GT's $\mu $ and $\nu $. Then the continuous functions
$f: (X, \mu ) \to (Y, \nu )$ are the following ones. For $|X| = 1 < |Y|$ there
is no continuous map $f: (X, \mu ) \to (Y, \nu )$. Else the continuous
maps are the (non-strictly) monotonically increasing or 
monotonically decreasing maps, which are continuous between 
the respective order topologies.}}


\begin{proof}
First we settle the case $\min \{ |X|, |Y| \} =0$. For $Y = \emptyset $ there
exists an $f: (X, \mu ) \to (Y, \nu )$ only if $X = \emptyset $. For 
$X = \emptyset $ we have $(X, \mu ) = (\emptyset, \{ \emptyset \} )$, hence to
any $(Y, \nu )$ there is exactly one set map 
$X \to Y$, that is continuous, and also is monotonous and continuous between
the respective order topologies.

For $\min \{ |X|, |Y| \} =1$ we may have
$|X| = |Y| = 1$ and then the unique set map $X \to Y$ is continuous, 
and also is monotonous and continuous between
the respective order topologies. 
For $|X| = 1 < |Y|$ we have that $(Y, \nu )$ is strong, while $X$
is not strong, hence there is no continuous map $f: (X, \mu ) \to (Y, \nu )$.
For $|X| > 1 = |Y|$ there is a unique set map $X \to Y$, that is
continuous, and also is monotonous and continuous between
the respective order topologies.

{\it{From now on let}} $|X|, |Y| \ge 2$.
Suppose that $f$ is not monotonous. Then there are $a,b,c,d \in X$ such that
\begin{equation}\label{equa4.1}
a<b {\text{ and }} f(a) < f(b), {\text{ and }} c<d {\text{ and }} f(c)>f(d)\,.
\end{equation}
Let $\{ a,b,c,d \} = \{
x_1, \dots , x_k \} $, where $x_i < x_{i+1}$. By (\ref{equa4.1}) 
$3 \le k \le 4$. If
for some $i$ we have
$f(x_i) = f(x_{i+1})$ then we delete from $ \{x_1, \dots , x_k \} $ the point
$x_{i+1}$. Thus we obtain $x'_1,\dots x'_l$, where by (\ref{equa4.1})
$x'_i < x'_{i+1}$ and  
$3 \le l \le 4$. Now for each $i$ we have 
$f(x_i) < f(x_{i+1})$ or $f(x_i) > f(x_{i+1})$.
By (\ref{equa4.1}) 
there is some $i$ such that $f(x_i) > f(x_{i+1}) < f(x_{i+2})$ or
$f(x_i) < f(x_{i+1}) > f(x_{i+2})$. We may suppose that we have the first case
(else we consider the converse ordering on $Y$). 

By Remark 4.14 the closed sets in $X$ and $Y$ are exactly the convex
subsets which are closed in the (finer)
order topology. Then the set $f^{-1} \left( \left[ \min \{ f(x_i),
\right. \right.$
\newline
$\left. \left. f(x_{i+2}) \} , \infty \right) \right) $ 
contains $x_i$ and $x_{i+2}$, but does not contain $x_{i+1}$. Thus the inverse
image of a closed set is not convex, hence is not closed, a contradiction.

So $f$ is monotonous. We may suppose that it is monotonically increasing (else
we take the converse ordering on $Y$). Then it suffices to prove that 
{\it{in the sense of topology, $f$ is
continuous from the left at any point}} $x \in X$. (Namely then a similar
reasoning shows that in the
sense of topology, $f$ is also continuous from the right, hence is, in the
sense of topology, continuous.) 
Observe that in the sense of topology $f$ is continuous from the left at $x
\in X$ if $x$ has an immediate predecessor in $X$. Therefore we suppose that 
$(-\infty , x)$ has no largest element (in particular, $x$ is not a minimal
element of $X$).

Suppose the contrary, namely that for some $x \in X$ and some $X \ni x' <x$, 
we have $f \left( (x',x) \right) \subset ( - \infty , y]$, for some $Y \ni
y<f(x)$, hence also $f \left( ( -
\infty , x) \right) \subset ( - \infty , y ]$, for some $Y \ni y < f(x)$. Then 
$(y, \infty )$ is a neighbourhood of $f(x)$ in the sense of GT's, such that
even no topological neighbourhood $N$ of $x$ satisfies $f\left( N \cap
(- \infty ,x) \right) \subset 
(y, \infty )$, hence no such $N$ satisfies $f(N) \subset (y, \infty )$,
although $(y, \infty )$ is also a topological neighbourhood of $f(x)$.

Conversely, continuity of a (non-strictly) monotonically increasing or
decreasing function
$f$ between the order topologies of $(X, \le )$ and $(Y, \le )$
implies its continuity between the order GT's 
of $(X, \le )$ and $(Y, \le )$ (recall $|X|, |Y| \ge 2$). 
It is sufficient to prove this
for $f$ monotonically increasing. It is sufficient to show that the inverse
image by $f$ of a basic open set in the order GT of $(Y, \le )$ 
is open in the order GT of $(X, \le )$. 
It is sufficient to prove this for a basic open set of the form 
$( - \infty ,y) \subset Y$. By continuity of $f$
in the topological sense we have
that $f^{-1}( - \infty ,y)$ is open in the topology of $(X, \le )$. Therefore
it is the union of non-empty basic open sets in the sense of topology, say
$f^{-1}( - \infty ,y) = \cup _{\alpha \in J} (a_{\alpha }, b_{\alpha
})$. Also, by the monotonically increasing property of $f$, we have that
$f^{-1}( - \infty ,y)$ is 
downward closed (i.e., $x' < x \in f^{-1}( - \infty
,y)$ implies $x' \in f^{-1}( - \infty ,y)$). Therefore we have also 
$f^{-1}( - \infty ,y) = \cup _{\alpha \in J} (- \infty , b_{\alpha })$. Thus 
$f^{-1}( - \infty ,y)$ is a union of open sets in the order GT of $(X, \le )$,
hence it is open in the order GT of $(X, \le )$ as well.
\end{proof}


For
topological spaces $X,Y$ where $Y$ is $T_2$, and continuous maps $f,g:X \to
Y$, if $f,g$ coincide on a dense subset of $X$, then they are equal. For GTS's
this is false.


\vskip.2cm

{\bf{Example 4.16.}} 
Let $X=Y=([0,1], \gamma _0)$. Let $f,g:[0,1] \to [0,1]$ be any
two different strictly
monotonically increasing maps in {\bf{GenTop}}, 
continuous in the topological sense, with
$f(0)=g(0)=0$ and $f(1)=g(1)=1$. Then $f,g$ are even homeomorphisms of
$([0,1], \gamma _0)$ in {\bf{GenTop}}, coinciding 
on the dense subset $\{0, 1 \} $ (cf. Remark 4.14), but $f \ne g$. 
Even, we can choose $f$ and $g$ so, that $\{ x \in [0,1] \mid f(x)=g(x) \} = \{
0,1 \} $.


\vskip.2cm

As well known (\cite{Eng}, Example 2.3.12 and Theorem 5.2.8, Hist. and
Bibl. Notes to \S 5.2) $T_4$ is not even finitely productive for
topological spaces. Also, all powers of some topological space $X$ 
are $T_4$ if and only if $X$ is compact $T_2$ (\cite{No}). 
However, for GTS's we have


\vskip.2cm

{\bf{Proposition 4.17.}} {\it{Let $\langle (X_{\alpha }, \mu _{\alpha }) \mid
\alpha \in J \} $ be an indexed set of normal (or $T_4$) GTS's. 
Then their product is also normal (or $T_4$).}}


\begin{proof}
Since $T_1$ is productive (cf. Proposition 4.7), 
we investigate the case of normal spaces only.

First we settle the case of the empty product. By the \S 2, {\bf{4}} this is the
GTS $(X_0, \{ \emptyset \} )$, where $|X_0|=1$, which is normal.

Now we suppose $J \ne \emptyset $.
Let $F_1,F_2$ be disjoint closed sets in $(X, \mu ) := \prod
_{\alpha \in J} (X_{\alpha }, \mu _{\alpha }) $. By Corollary 3.10
we have $F_i=
\prod _{\alpha \in J} F_{i, \alpha }$, where $F_{i, \alpha } \subset 
X_{\alpha }$ is $\mu _{\alpha }$-closed.
If for each $\alpha \in J$ we have $F_{1, \alpha } \cap F_{2, \alpha } \ne
\emptyset $, then $F_1 \cap F_2 \ne \emptyset $. Hence for some $\alpha _0 \in
J$ we have $F_{1, \alpha _0} \cap F_{2, \alpha _0} = \emptyset $. Then by
normality of $(X_{\alpha _0}, \mu _{\alpha _0})$ there are disjoint
$\mu _{\alpha _0}$-open sets $M_{1, \alpha _0}, M_{2, \alpha _0}$, such that
$F_{i, \alpha _0} \subset M_{i, \alpha _0}$ for $i=1,2$. Then $F_i \subset
M_{i, \alpha _0} \times \prod _{\alpha \in J \setminus \{ \alpha _0 \} } X_
{\alpha }$, for $i=1,2$. These last sets are disjoint $\mu $-open sets in $X$. 
\end{proof}


We say that a GTS $(X, \mu )$ is {\it{compact}}, or {\it{Lindel\"of}}
if any open cover of $X$ has a
finite, or at most countably infinite subcover of $X$, respectively. This
is the exact analogue of compactness and the Lindel\"of property for
topological spaces. In particular, if $X \not\in \mu $, then $X$ is compact.
More generally, for $\kappa $ an infinite cardinal, we say that 
a GTS $(X, \mu )$ is $\kappa $-{\it{compact}} if any open cover of 
$(X, \mu )$ has an (open) 
subcover of $X$, of cardinality less than $\kappa $. This is the
analogue of $\kappa $-compactness for topological spaces, cf. \cite{Ju},
p.6. For $\kappa = \aleph _0$ this is compactness, for $\kappa =\aleph _1$
this is the Lindel\"of property.


By Tychonoff's theorem compactness is {\it{productive}} 
for topological spaces (i.e., products of compact spaces are compact), but
Lindel\"ofness is not, cf. \cite{Eng}, 3.8.15. For GTS's the situation is
very different, as shown by the next Proposition.


\vskip.2cm

{\bf{Proposition 4.18.}} {\it{For any infinite cardinal $\kappa $, $\kappa
    $-compactness is closed-he\-re\-di\-ta\-ry and productive for
GTS's, and is also inherited by surjective images for GTS's.}}


\begin{proof}
Closed-hereditariness and inheriting by surjective images 
are proved exactly as for topological spaces.
 
We turn to products. Let $\langle (X_{\alpha } , \mu _{\alpha }) 
\mid \alpha \in J \rangle $ be an indexed set of $\kappa $-compact 
GTS's, with product $(X, \mu )$. We may suppose that each $ X_{\alpha }$ is
non-empty. 
Let ${\cal{G}} =
\{ G_{\beta } \mid \beta \in B\} $ be an open cover of $(X, \mu )$.
By Corollary 3.10 we have for each $\beta \in B$ that
$$
G_{\beta } = \cup _{\alpha \in J} \pi _{\alpha } ^{-1} (M_{\alpha , \beta }) \,,
$$
where $M_{\alpha , \beta } \in \nu _{\alpha }$. Also by Corollary 3.10
we have that
an open base of $(X, \mu )$ is $\cup \{ \pi _{\alpha }^{-1}(\nu _{\alpha
}) \mid \alpha \in J \} $. 
We define the set system ${\cal{H}}$ on $X$ as follows:
$$
{\cal{H}} := \{ \pi _{\alpha } ^{-1} (M_{\alpha , \beta }) \mid \alpha \in
J,\,\,\,\,\beta \in B \} \,.
$$
Then ${\cal{H}}$ 
is a cover of $X$ since $\cup {\cal{H}} = \cup {\cal{G}}$ and 
${\cal{G}}$ is a cover of $X$, and consists of open sets in $(X, \mu )$ since
each element of ${\cal{H}}$ belongs to the above mentioned base of $(X, \mu
)$. 

We distinguish two cases. 
\newline
(i)
Either for each $\alpha \in J$ we have $\cup \{ 
M_{\alpha , \beta } \mid \beta \in B \} \ne X_{\alpha }$, or
\newline
(ii) 
for some $\alpha _0 \in J$ we have 
$\cup \{ M_{\alpha _0, \beta } \mid \beta \in B \} = X_{\alpha _0}$.
\newline
In case (i) we choose for each $\alpha \in J$ a point 
$x_{\alpha } \in X_{\alpha }$ not
contained by the union in (i). Then the point $\langle x_{\alpha } \mid \alpha
\in J \rangle $ is not covered by ${\cal{H}}$, a contradiction. 
In case
(ii) the union in (ii) is an open cover of $X_{\alpha _0}$, hence by $\kappa
$-compactness of $X_{\alpha _0}$ it has an open
subcover $ \{ M_{\alpha _0, \beta } \mid \beta \in B' \} $ 
of $X_{\alpha _0}$ with $|B'| < \kappa $. Then  
$\{ \pi _{\alpha  _0}^{-1} (M_{\alpha _0, \beta }) \mid \beta \in B' \} $
is an open cover of $(X, \mu )$, of cardinality less than $\kappa $.
Now recall that for each $\beta \in B' \subset B$ we have 
$\pi _{\alpha  _0}^{-1} (M_{\alpha _0, \beta }) \subset G_{\beta }$. Therefore
$\{ G_{\beta } \mid \beta \in B' \} $ is a subset of
${\cal{G}}$, which is also an open cover of $X$, and has cardinality
less than $\kappa $.
\end{proof}


A topological space $X$ is $T_{3.5}$ if and only if it is a subspace or a
dense subspace of some (compact)
$T_4$ space. Namely, one can consider the Stone-\v Cech compactification of
$X$ that is compact $T_2$ hence $T_4$. It is interesting that an analogue of
this statement holds also for GTS's, although with a completely
different proof.


\vskip.2cm

{\bf{Proposition 4.19.}}
{\it{A GTS $(X, \mu )$ is $T_{3.5}$ (completely regular) 
if and only if it is homeomorphic to a subspace, or
equivalently to a dense subspace of a $T_4$ GTS, which can be supposed
to be also compact  (has the weak
structure w.r.t. a map, or
equivalently w.r.t. a dense map from $UX$ to a $T_4$ GTS,
which can be supposed to be also compact).}}


\begin{proof}
{\bf{1.}}
We begin with the ``if'' part.
$T_4$ implies the hereditary property $T_{3.5}$ (cf. Proposition 4.7), 
which implies complete regularity,
and complete regularity is inherited by initial (weak) structures, in
particular for sources consisting of one map, cf. Proposition 4.7. 
Hence all subspaces of $T_4$ spaces are $T_{3.5}$, and initial (weak) 
structures for all sources
consisting of a single map $\varphi :X \to U(Y, \mu )=Y$ with $Y$ $T_4$ are
completely regular.

{\bf{2.}}
Conversely, for the ``only if'' part,
we begin with the trivial cases $|X| \le 1$, i.e., with $(
\emptyset , \{ \emptyset \})$, and with $(X, \{ \emptyset \})$ and 
$(X, \{ \emptyset , X \} )$ where $|X|=1$. Each of these three spaces are finite
hence compact, are $T_1$ and normal (the second space has no disjoint closed
sets, and the other ones are discrete, i.e., of the form $\left( X, P(X)
\right) $, which are normal). Then they are dense subspaces of themselves,
which proves for them the ``only if'' part of the theorem.

Now let $(X, \mu )$ be $T_{3.5}$ with $|X| \ge 2$.
By Theorem 4.12, $(X, \mu )$ is homeomorphic 
to a dense subspace of some power of $([0,1], \gamma _0)$.
By Remark 4.1, $([0,1], \gamma _0)$ is $T_4$, and then by Proposition 4.17
each power of $([0,1], \gamma _0)$ is $T_4$. Since
$([0,1], \gamma _0)$ is coarser than the topological $[0,1]$ space,
it is compact as well, hence all its powers are compact as well by Proposition
4.18. This ends the proof of the ``only if'' part for $T_{3.5}$ GTS's.

{\bf{3.}}
We turn to complete regularity, to the ``only if'' part.
Let $(X, \mu )$ be completely regular, let $(Y, \nu )$ be its
$T_0$-reflection, with canonical quotient map $q:(X, \mu ) \to (Y, \nu )$ from
Definition 4.10. Then by Proposition 4.11, (2) $(Y, \nu )$ is $T_{3.5}$,
hence by {\bf{2}} of this proof it admits a dense embedding $i:(Y,
\nu ) \to (Z, \varrho )$ to some compact $T_4$ GTS $(Z, \varrho )$. Then the
source consisting of the single map $U(iq) : U(X, \mu ) = X \to U(Z, \varrho )$
has as initial lift (weak structure) $(X, \mu )$, and $iq:(X, \mu ) \to
(Z, \varrho )$ is a dense map.
This ends the proof of the ``only if'' part for completely regular GTS's.
\end{proof}


As known, a $T_2$ topological space of density $\kappa $ has cardinality at
most ${\text{exp}}\left( {\text{exp}} ( \kappa ) \right) $, cf. \cite{Ju},
p. 13. For GTS's
the situation is completely different.

\vskip.2cm

{\bf{Example 4.20.}}
For GTS's we may have a two-point dense subset in a compact $T_4$ 
GTS of arbitrarily large
cardinality. Namely, in $([0,1], \gamma _0) ^\kappa $, which is compact by 
Proposition 4.18 and is $T_4$ by Proposition 4.17,
the two points having
all coordinates $0$, or all coordinates $1$, form a dense subspace in 
$([0,1], \gamma _0) ^\kappa $ (cf. Corollary 3.10).


\section{Subspaces and sums of GTS's}


We recall Corollary 3.13
about the construction of the sum of an indexed set $\langle (X
_{\alpha }, \mu _{\alpha}) \mid \alpha \in J \rangle $ 
or $\langle (X_{\alpha },c_{\alpha }) \mid \alpha \in J \rangle $ of GTS's.

As mentioned in \S 2, {\bf{5}}, a
common way to produce GT's is given by the GT's
$\mu ( \gamma )$ (see \cite{Cs0.5} and also
our \S 2, {\bf{5}}), where $X$ is a set
and $\gamma \in \Gamma (X)$.
 
Let $\gamma _{\alpha } :X_{\alpha } \to X_{\alpha }$ be monotonous. We define
$\gamma : P \left( \coprod _{\alpha \in J} X_{\alpha } \right) \to P \left(
\coprod _{\alpha \in J} X_{\alpha } \right) $ 
by 
$$
\gamma (A) := \coprod _{\alpha \in J}
\gamma _{\alpha } (A \cap X_{\alpha })\,.
$$
Then $\gamma \in \Gamma \left( \coprod _{\alpha \in J} X_{\alpha } \right) $.

Now we can consider two GT's on the same set $X := \coprod_{\alpha\in
J}X_{\alpha}$. The first one is $\mu (\gamma )$ and the second one is
$\coprod _{\alpha \in J}  \mu ( \gamma _{\alpha })$ 
(notations cf. in \S 2, {\bf{5}}). 

Similarly, for $\gamma : P(X) \to P(X)$ monotonous, and $X_0 \subset X$, for
$\gamma _0 : P(X_0) \to  P(X_0)$ defined by 
$$
\gamma _0 A_0 := (\gamma A_0) \cap X_0 \,,
$$
we have $\gamma _0 \in \Gamma (X_0)$. Then we have two GT's on
$X_0$, the first one is $\mu (\gamma _0)$ and the second one is $\mu (\gamma
)|X_0$.

We are going to
compare these two pairs of GT's by the following proposition.
We recall that $\gamma : P(X) \to P(X)$ is {\it{completely additive}} if for
any $\{ A_{\alpha } \mid \alpha \in J \} \subset P(X)$ we have $\gamma \left(
\cup _{\alpha \in J} A_{\alpha } \right) = 
\cup _{\alpha \in J} \gamma (A_{\alpha })$.


\vskip.2cm

{\bf{Proposition 5.1.}}
{\it{(1) Let $\langle X_{\alpha } \mid \alpha \in J \rangle $ be
an indexed set of sets
and $\gamma _{\alpha }: P(X_{\alpha }) \to P(X_{\alpha })$, 
for $\alpha \in J$ be monotonous.
Then, with the notations introduced before this Proposition,
$$
\left( \coprod _{\alpha \in J} X_{\alpha }, \mu (\gamma ) \right)
=
\coprod _{\alpha \in J} \left( X_{\alpha }, \mu (\gamma _{\alpha }) \right) \,.
$$
(2) Let $X$ be a set and $\gamma : P(X) \to P(X)$ be monotonous, and 
let $X_0 \subset X$.
Then, with the notations introduced before this Proposition, 
$$
\mu (\gamma _0) \subset \mu (\gamma ) |X_0 \,.
$$ 
The converse inclusion is false even if $\gamma $ is completely additive and
$\gamma \left( P(X) \right) = \{ \emptyset , X \} $.
}}


\begin{proof}
{\bf{1.}}
As in \S 2, {\bf{4}}, we suppose that 
the sets $X_{\alpha }$, for $\alpha \in J$,
form a partition of $\coprod _{\alpha \in J} X_{\alpha }$.
Let $A \subset \coprod _{\alpha \in J}$. We write $A_{\alpha } := A \cap
X_{\alpha }$; hence $A=\coprod _{\alpha \in J} A_{\alpha }$. We have $A
\subset \gamma A \Longleftrightarrow \forall \alpha \in J \,\,\,\,A_{\alpha }
\subset \gamma _{\alpha }A_{\alpha }$. Hence $\mu (\gamma ) = \coprod _{\alpha
\in J} \mu ( \gamma _{\alpha })$.

{\bf{2.}}
First we show $\mu ( \gamma _0)  \subset \mu ( \gamma )|X_0$. We have for $A_0
\subset X_0$ that $A_0 \in \mu ( \gamma _0) \Longleftrightarrow 
A_0 \subset \gamma _0 A_0 \Longleftrightarrow A_0 \subset \gamma (A_0) \cap
X_0 \Longleftrightarrow A_0 \subset \gamma (A_0) \Longleftrightarrow
A_0 \in \mu ( \gamma )$. Then $A_0 \in \mu ( \gamma _0) \Longleftrightarrow
A_0 \in \mu ( \gamma ) \Longrightarrow A_0 = A_0 \cap X_0 \in 
\mu ( \gamma )|X_0$, therefore $\mu (\gamma _0 ) \subset \mu ( \gamma )|X_0$.

About the converse inclusion we give the following counterexample. Let $|X|
\ge 2$, $X_0 \subset X$ and $X_0 \not\in \{ \emptyset , X \} $. Let $\gamma
\in \Gamma (X)$ be defined as follows. For $A \subset X_0$ we have $\gamma
(A)  = \emptyset $, for $A \subset X$ and $A \not\subset X_0$ we have $\gamma
(A) = X$. Thus $\gamma $ is completely additive and $\gamma \left( P(X)
\right) = \{ \emptyset , X \} $. Then by $X_0 \ne X$ we have
$\mu ( \gamma ) = \{ A \subset X \mid A \subset \gamma A \} = \{ \emptyset \}
\cup \{ A \subset X \mid A \not\subset X_0 \} $, hence 
$$
\mu ( \gamma ) |X_0 = P(X_0)\,.
$$
On the other hand
\[
\begin{array}{c}
\mu (\gamma _0) = \{ A_0 \subset X_0 \mid A_0 \subset \gamma _0 A_0 \}  =  \{
A_0 \subset X_0 \mid A_0 \subset \gamma (A_0) \cap X_0 \} = \\
\{ A_0 \subset
X_0 \mid A_0 \subset \gamma (A_0) \} =  \{ \emptyset \} 
\end{array}
\]
Therefore, by $X_0 \ne \emptyset $ we have
$$
\mu ( \gamma _0)= \{ \emptyset \} 
\not\supset P(X_0)= \mu (\gamma )|X_0\,.
$$
\end{proof}


As mentioned in \S 2, {\bf{5}}, another
common way to produce GT's is given by the GT's
$\kappa (\mu , k )$ (see \cite{Cs3} and also our \S 2, {\bf{5}}), 
where $(X, \mu )$ is a GTS and 
$k: \mu \to P(X)$ is an enlargement on $(X, \mu )$.

Now, let $\langle (X_{\alpha},\mu_{\alpha}) \mid \alpha\in J \rangle $ 
be an indexed set of 
GTS's and $k_{\alpha}:\mu_{\alpha}\rightarrow P(X_{\alpha})$ be an
enlargement on $X_{\alpha}$, for $\alpha\in J$. Let
\begin{equation}
\begin{array}{c}
(X, \mu ) := \coprod _{\alpha \in J} (X_{\alpha }, \mu _{\alpha }) 
{\text{ and}} \\
\left( \coprod _{\alpha \in J} X_{\alpha }, \coprod _{\alpha \in J} \kappa 
( \mu _{\alpha }, k_{\alpha }) \right) := \coprod _{\alpha \in J} \left(
X_{\alpha }, \kappa ( \mu _{\alpha }, k_{\alpha } ) \right) 
\end{array}
\end{equation}
We can ask if there is
an enlargement $k$ on the sum set $X := \coprod_{\alpha\in
J}X_{\alpha}$, such that $\left( X, \kappa (\mu, k) \right) =
\coprod _{\alpha \in J}
\left( X_{\alpha}, \kappa (\mu_{\alpha}, k_{\alpha }) \right)$.
Actually, here we will have only a one-sided inclusion. The converse
inclusion is in general false, but the necessary and sufficient
condition for equality will be given.

Similarly, for subspaces $(X_0, \mu _0)$ of $(X, \mu )$, one could ask if
on $(X_0, \mu _0)$ there is an enlargement $k_0 : \mu _0 \to P(X_0)$, such that 
$\left( X_0, \kappa ( \mu _0, k_0) \right) = \left (X, \kappa ( \mu , k )
\right) |X_0$. Actually, here we will have only a one-sided inclusion,
and only under some additional hypotheses. The converse inclusion is false, even
under more restrictive additional hypotheses.


\vskip.2cm

{\bf{Proposition 5.2.}}
{\it{Let $\langle (X_{\alpha},\mu_{\alpha}) \mid {\alpha\in J} \rangle $ 
be an indexed set of 
GTS's and $k_{\alpha} : \mu_{\alpha} \rightarrow P(X_{\alpha})$ be an
enlargement on $(X_{\alpha}, \mu _{\alpha })$, for $\alpha\in J$. 
Then the enlargement
$k$ on $(X, \mu ) := \coprod_{\alpha\in J}(X_{\alpha },\mu_{\alpha})$ 
defined for $M \in \mu $ 
by $kM := \cup _{\alpha \in J} k_{\alpha } (M \cap X_{\alpha })$
satisfies 
$\kappa ( \mu ,k) \subset \coprod _{\alpha \in J} \kappa ( \mu _{\alpha },
k_{\alpha })$. 
Here equality 
holds if and only if either 
$| \{ \alpha \in J \mid X_{\alpha } \ne \emptyset \} | \le 1$ or 
$| \{ \alpha \in J \mid X_{\alpha } \ne \emptyset \} | \ge 2$ and for each
$\alpha \in J$ we have $k_{\alpha } \emptyset = \emptyset $.

Let $(X, \mu )$ be a GTS and $k :\mu \to P(X)$ be an enlargement on $(X, \mu )$.
Let $X_0 \subset X$.
If $k$ is also monotonous and $\mu $ is closed under the intersections of
pairs of elements and $X_0$ is $\mu $-open, then 
the enlargement $k_0$ on $(X_0, \mu _0)$ defined for $M_0 \in \mu _0$
by $k_0M_0:= (kM_0)
\cap X_0$ satisfies $ \kappa ( \mu , k) | X_0 \subset
\kappa ( \mu _0, k_0) $. The converse inclusion is false
even for $\mu $ a topology, $X_0$ $\mu $-open, and $k$ a topological closure.}}


\begin{proof}
{\bf{1.}}
We begin with the case of sums.
As in \S 2, {\bf{4}}, we suppose that the sets $X_{\alpha }$ for 
$\alpha \in J$ form a partition of $\coprod _{\alpha \in J} X_{\alpha }$.

We have
$kM = \cup_{\alpha\in J} k_{\alpha} (M_{\alpha}) $, where $M \in \mu $ and
$M_{\alpha } := M \cap X_{\alpha }$. Then $M = \coprod _{\alpha \in J}
M_{\alpha } \subset \coprod _{\alpha \in J} k_{\alpha } M_{\alpha } = kM
\subset X$, thus $k$ is an enlargement on $(X, \mu )$.

Since for $X_{\alpha } = \emptyset $ we have $\mu _{\alpha } = \kappa (\mu
_{\alpha }, k_{\alpha }) = \{ \emptyset \} $, these contribute nothing to 
either $\mu $, or to $k$, or to $\kappa ( \mu , k)$, or to 
$\coprod _{\alpha \in J} \kappa (
\mu _{\alpha } , k_{\alpha })$, therefore we may omit all empty $X_{\alpha }$'s
simultaneously. Therefore we will suppose that each $X_{\alpha }$ is non-empty.
If we have at most one non-empty $X_{\alpha }$, then evidently $\coprod
_{\alpha \in J} \kappa ( \mu _{\alpha }, k_{\alpha }) = \kappa (\mu , k)$.
Therefore we suppose that there are at least two non-empty $X_{\alpha }$'s.

First we prove 
\begin{equation}\label{equaalpha}
\kappa (\mu ,k ) \subset 
\coprod _{\alpha \in J} \kappa (\mu _{\alpha },k_{\alpha } )\,.
\end{equation}
Let $A \in \kappa ( \mu , k) \subset \mu $. We write $A_{\alpha }
:=  A \cap X_{\alpha } \in \mu _{\alpha }$. If $A_{\alpha } = \emptyset $, 
then $A_{\alpha } \in \kappa ( \mu _{\alpha }, k_{\alpha })$, so 
we need to deal only with such $\alpha $'s, for
which $A_{\alpha } \ne \emptyset $. Let $x_{\alpha } \in A_{\alpha } \subset
A$. Then there exists an $M \in \mu$ such that
$x_{\alpha } \in M$ and $kM \subset
A$. Writing $M_{\alpha } := M \cap A_{\alpha }$ and $M_{\beta } := M \cap
X_{\beta }$ for $\beta \in J
\setminus \{ \alpha \} $, this means that
$x_{\alpha } \in M_{\alpha } \in \mu _{\alpha }$ and $k_{\alpha } M_{\alpha }
\subset A_{\alpha }$ and for $\beta \in J \setminus \{ \alpha \} $ that
$M_{\beta } \in \mu _{\beta }$ and
$k_{\beta } M_{\beta } \subset A_{\beta }$.
(For $\beta \ne \alpha $ the
condition for $A_{\beta } $ is weaker than the condition for $A_{\alpha }$, 
so we disregard the condition about all
$\beta \in J \setminus \{ \alpha \} $. Observe that any $\beta \in J \setminus
\{ \alpha \} $ --- with $X_{\beta } \ne \emptyset $ ---
can occur in the role of $\alpha $.)
Then we obtain $A_{\alpha } \in
\kappa ( \mu _{\alpha }, k_{\alpha })$, that is, (\ref{equaalpha}) is proved.

Second we deal with the validity of the inclusion 
$$
\coprod _{\alpha \in J} \kappa ( \mu _{\alpha } , k_{\alpha }) \subset 
\kappa ( \mu , k) \,.
$$
Since $\coprod _{\alpha \in J} \kappa ( \mu _{\alpha },
k_{\alpha })$ has as base $\cup _{\alpha \in J} \kappa ( \mu _{\alpha },
k_{\alpha })$, this inclusion is equivalent to 
$$
\cup _{\alpha \in J} \kappa ( \mu _{\alpha }, k_{\alpha }) \subset \kappa 
(\mu , k), {\text{ i.e., to }} \forall \alpha \in J\,\,\,\,
\kappa ( \mu _{\alpha }, k_{\alpha }) \subset \kappa (\mu , k)\,. 
$$
Let
$A_{\alpha } \in \kappa ( \mu _{\alpha }, k_{\alpha })$, i.e., 
\begin{equation}\label{equagamma}
\begin{array}{l}
{\text{for each }}
x_{\alpha } \in A_{\alpha } {\text{ there exists an }}
M_{\alpha } \in \mu _{\alpha } \\
{\text{such that }}
x_{\alpha } \in M_{\alpha } \subset
k_{\alpha } M_{\alpha } \subset A_{\alpha } \subset X_{\alpha }\,. 
\end{array}
\end{equation}
Then $A_{\alpha } \in \kappa ( \mu , k)$ if and only if 
\begin{equation}\label{equadelta}
\begin{array}{l}
{\text{for each }}
x_{\alpha } \in A_{\alpha } {\text{ there exists an }} M \in \mu \\
{\text{such that }} x_{\alpha } \in M \subset kM \subset A_{\alpha } \subset
X_{\alpha }\,. 
\end{array}
\end{equation}
Here $M \in \mu $ and $M \subset X_{\alpha }$ means $M =:
M^{\alpha} \in \mu _{\alpha }$, hence (\ref{equadelta}) is equivalent to
\begin{equation}\label{equai.j}
\begin{array}{l}
{\text{for each }}
x_{\alpha } \in A_{\alpha } {\text{ there exists an }} M \in \mu 
{\text{ such that}} \\
X_{\alpha } \supset A_{\alpha } \supset
kM = kM^{\alpha } = 
\cup _{\beta \in J} k_{\beta } M^{\alpha } \\
= (k_{\alpha }M^{\alpha }) \coprod
\left( \coprod _{\beta \in J \setminus \{ \alpha \} } k_{\beta } \emptyset
\right) \supset M^{\alpha } = M \ni x_{\alpha }
\,.
\end{array}
\end{equation} 
Then $k_{\alpha }M_{\alpha } \subset A_{\alpha }$ 
is satisfied for $M^{\alpha } :=
M_{\alpha } $, thus we need to satisfy yet that for all 
$\beta \in J \setminus \{ \alpha \} $ we have for $k_{\beta } \emptyset
\subset X_{\beta }$ that also $k_{\beta } \emptyset \subset X_{\alpha }$ (cf. 
(\ref{equai.j}).
But then $k_{\beta } \emptyset \subset X_{\alpha } \cap X_{\beta }
=  \emptyset $, i.e., $k_{\beta } \emptyset = \emptyset $.
However, this is just the necessary and
sufficient hypothesis for the case $| \{ \alpha \in J \mid X_{\alpha } \ne
\emptyset \} | \ge 2$,
given in this Theorem. (Observe that here we have
$k_{\beta } \emptyset =  \emptyset $ only for $\beta \in J \setminus \{
\alpha \} $. However, since now $|J| \ge 2$, we change $\alpha \in J$ to another
element $\alpha ' \in J$, and
then we obtain $k_{\alpha } \emptyset =  \emptyset $ as well.)

Conversely, $| \{ \alpha \in J \mid X_{\alpha } \ne \emptyset \} | \ge 2$ and
$\forall \alpha \in J\,\,\,\,k_{\alpha } \emptyset = \emptyset $ 
implies $\kappa ( \mu _{\alpha }, k_{\alpha } ) $
\newline
$\subset \kappa (\mu ,k)$, 
i.e., (\ref{equagamma}) 
$\Longrightarrow $ (\ref{equai.j}) $\Longleftrightarrow $ (\ref{equadelta}).

{\bf{2.}}
We turn to the case of subspaces.
We have $k_0M_0 =  k (M_0) \cap X_0$
for $M_0 \in \mu _0$ (thus $M_0 \subset X_0$). Here $k (M_0)$ is defined since
$M_0$ is the intersection of some open set of $(X, \mu )$ and $X_0$, hence is 
$\mu $-open by hypothesis.
Then $k_0 M_0 \supset  M_0 \cap X_0 = M_0$,
thus $k_0$ is an enlargement on $(X_0, \mu _0)$.

Let $A \subset X$ be $\kappa (\mu , k )$-open, i.e., $x \in A$ implies
$\exists M \in \mu $ such that
$x \in M$ and $kM \subset A$. Let $A_0:=A \cap
X_0$. Then for $x_0 \in A_0 \,\,\,\,(\subset X_0)$ there exists 
$x_0 \in M \subset
A $ such that $kM \subset A$. Then also $x_0 \in M \cap X_0 \subset A \cap X_0$
and $k(M \cap X_0) \subset kM \subset A$ by monotony of $k$. 
Then also $k_0(M \cap X_0) =
k(M \cap X_0) \cap X_0 \subset A \cap X_0$ while $M \cap X_0 \in \mu
_0$ by hypothesis. Hence $A_0 = A \cap X_0$ is $\kappa (\mu _0, k_0)$-open.

About the converse inclusion we give the following counterexample. 

Let $X := \{ 1/n \mid n \in {\mathbb{N}} \} \cup \{ 0 \} $ and $X_0 := 
\{ 1/n \mid n \in {\mathbb{N}} \} $ its open subspace, with the usual
topologies. We define $k : P(X) \to P(X)$ as follows: $k \emptyset = \emptyset
$, and for $\emptyset \ne M \subset X$ we define $kM := M \cup \{ 0 \} $. This
is a topological closure. 

Then for $A \subset X$ we have $A \in \kappa ( \mu , k)$ if and only if 
\begin{equation}\label{equa5.1}
\begin{array}{l}
x \in A \Longrightarrow \exists M \in \mu  {\text{ such that }} x
\in M {\text{ and }} kM = M \cup \{ 0 \} \subset A \,.
\end{array}
\end{equation}
Now we will determine $\kappa ( \mu , k)$. 
\newline
(1) If $A = \emptyset $ then $A \in \kappa ( \mu , k)$.
\newline
(2) If $A \ne \emptyset $ and $0 \not\in A$ then supposing (\ref{equa5.1})
we have $\exists x \in A$, and then 
$\{ 0 \} \subset A$, which is a contradiction. Hence such an $A$ does not
belong to $\kappa ( \mu , k)$.
\newline
(3) If $A \ne \emptyset $ and $0 \in A$, then $x \in A$ can be either $1/n$
for some $n \in {\mathbb{N}}$, or it can be $0$. For $x=1/n$ (\ref{equa5.1})
is satisfied for $M := \{ 1/n \} $, since $0 \in A$. For $x=0$ (\ref{equa5.1})
means $\exists m \in {\mathbb{N}}$ such that $0$ has a neighbourhood 
$0 \in \{ 1/m, 1/(m+1), \ldots \} \cup \{  0  \} $ in $X$, such that
$\{ 1/m, 1/(m+1), \ldots \} \cup \{  0  \} \subset A$.
This means that $0 \in A$ and $A \setminus \{ 0 \} $ is cofinite in $X
\setminus \{ 0 \} = X_0$.
That is, 
\begin{equation}\label{equa5.2}
\begin{array}{l}
\kappa ( \mu , k) = \{ \emptyset \} \cup \{ A \subset X \mid 0 \in A {\text{
and}} \\
A \setminus \{ 0 \} {\text{ is cofinite in }} X
\setminus \{ 0 \} = X_0 \} \,.
\end{array}
\end{equation}
Hence 
\begin{equation}\label{equa5.3}
\begin{array}{l}
\kappa ( \mu , k)|X_0 {\text{ is the cofinite topology on }} X_0 \,.
\end{array}
\end{equation}

Turning to $X_0$, we have that $\mu _0$ is the discrete topology on $X_0$,
and, for $M_0 \subset X_0$, we have for $M_0 = \emptyset $ that
$k_0 \emptyset := k( \emptyset ) \cap X_0 = \emptyset $, and for
$M_0 \ne \emptyset $ we have
$k_0M_0 := k(M_0) \cap X_0 = (M_0 \cup \{ 0 \} ) \cap X_0 =
M_0$. Therefore $k_0$ is the closure associated to the discrete topology on
$X_0$. Hence each $A_0 \subset X_0$ is $\kappa ( \mu _0, k_0)$-open, since
for $x_0 \in A_0$ we can choose $M_0 := \{ x_0 \} $ and then $x_0 \in 
\{ x_0 \} $ and $k_0 \{ x_0 \} = \{ x_0 \} \subset A_0$. That is, we have
\begin{equation}\label{equa5.4}
\begin{array}{l}
\kappa (\mu _0, k_0) = P(X_0) \,.
\end{array}
\end{equation}

By (\ref{equa5.3}) and (\ref{equa5.4}) we have
\begin{equation}\label{equa5.5}
\begin{array}{l}
\kappa (\mu _0, k_0) \not\subset \kappa ( \mu , k)|X_0 \,.
\end{array}
\end{equation}
\end{proof}


{\bf{Remark 5.3.}}
The construction of the counterexample is a special case of the $\theta
$-modification of a bitopological space (which is itself a special case of the
$\theta $-modification of a bi-GTS, cf. \cite{MAH} and \cite{CsM}).



\noindent
E. Makai, Jr.\\
MTA Alfr\'ed R\'enyi Institute of Mathematics\\
H-1364 Budapest, Pf. 127, Hungary\\
{\rm{http://www.renyi.mta.hu/\~{}makai}}\\ 

\noindent
Esmaeil Peyghan and Babak Samadi\\
Faculty  of Science, Department of Mathematics\\
Arak University\\
Arak 38156-8-8349,  Iran\\

\noindent
Email: makai.endre@renyi.mta.hu,\\
e-peyghan@araku.ac.ir, \ \ b\_samadi61@yahoo.com

\end{document}